\documentclass{amsart}
\usepackage[all]{xy}
\usepackage{verbatim}
\usepackage{color}
\usepackage{amsthm}
\usepackage{amssymb}
\usepackage[colorlinks=true]{hyperref}
\usepackage{graphicx}
\usepackage{mathtools}
\usepackage[utf8]{inputenc}

\usepackage{footmisc}
\usepackage{stmaryrd}

\numberwithin{equation}{section}

\newtheorem{theorem}[equation]{Theorem}
\newtheorem*{theorem*}{Theorem}

\newtheorem*{conjecture*}{Mamma Conjecture}
\newtheorem*{conjecture1*}{Mamma Conjecture (revisited)}

\newtheorem{corollary}[equation]{Corollary}
\newtheorem*{corollary*}{Corollary}

\theoremstyle{remark}

\newtheorem{example}[equation]{Example}

\newtheorem{notation}[equation]{Notation}

\theoremstyle{remark}
\newtheorem{remark}[equation]{Remark}

\newcommand{\cA}{{\mathcal A}}
\newcommand{\cB}{{\mathcal B}}
\newcommand{\cC}{{\mathcal C}}
\newcommand{\cD}{{\mathcal D}}

\newcommand{\cF}{{\mathcal F}}

\newcommand{\cL}{{\mathcal L}}
\newcommand{\cM}{{\mathcal M}}

\newcommand{\cO}{{\mathcal O}}

\newcommand{\cT}{{\mathcal T}}

\newcommand{\cW}{{\mathcal W}}
\newcommand{\cX}{{\mathcal X}}
\newcommand{\cY}{{\mathcal Y}}
\newcommand{\cZ}{{\mathcal Z}}

\newcommand{\bbA}{\mathbb{A}}

\newcommand{\bbC}{\mathbb{C}}

\newcommand{\bbF}{\mathbb{F}}
\newcommand{\bbG}{\mathbb{G}}

\newcommand{\bbP}{\mathbb{P}}

\newcommand{\bbQ}{\mathbb{Q}}
\newcommand{\bbZ}{\mathbb{Z}}

\DeclareMathOperator{\id}{id}

\DeclareMathOperator{\NChow}{NChow} % category of noncommutative Chow motives
\newcommand{\dgcat}{\mathrm{dgcat}} % codimension 
\newcommand{\perf}{\mathrm{perf}}

\newcommand{\Chow}{\mathrm{Chow}}

\newcommand{\dg}{\mathrm{dg}}

\newcommand{\Hom}{\mathrm{Hom}}

\newcommand{\dgHo}{\mathrm{H}^0}

\newcommand{\op}{\mathrm{op}}

\newcommand{\too}{\longrightarrow}

\newcommand{\ie}{\textsl{i.e.}\ }

\let\oldmarginpar\marginpar
\def\marginpar#1{\oldmarginpar{\tiny #1}}

\def\multiset#1#2{\ensuremath{\left(\kern-.3em\left(\genfrac{}{}{0pt}{}{#1}{#2}\right)\kern-.3em\right)}}

\begin{document}

\title[Noncommutative counterparts of celebrated conjectures]{Noncommutative counterparts \\of celebrated conjectures}
\author{Gon{\c c}alo~Tabuada}

\address{Gon{\c c}alo Tabuada, Department of Mathematics, MIT, Cambridge, MA 02139, USA}
\email{tabuada@math.mit.edu}
\urladdr{http://math.mit.edu/~tabuada}
\thanks{The author is very grateful to Guillermo Corti\~nas for the invitation to participate in the conference {\em $K$-theory in algebra, analysis and topology}, Buenos Aires, Argentina. The author was supported by a NSF CAREER Award.}

\date{\today}

\abstract{In this survey, written for the proceedings of the conference {\em $K$-theory in algebra, analysis and topology}, Buenos Aires, Argentina (satellite event of the ICM 2018), we give a rigorous overview of the noncommutative counterparts of some celebrated conjectures of Grothendieck, Voevodsky, Beilinson, Weil, Tate, Parshin, Kimura, Schur, and others.}}

\maketitle

%\smallskip

%\vskip-\baselineskip
%\vskip-\baselineskip
%\vskip-\baselineskip
%\tableofcontents

%\vspace{-0.3cm}

\tableofcontents

%-------------------------------------------------------------------------------
\section*{Introduction}
%-------------------------------------------------------------------------------
Some celebrated conjectures of Grothendieck, Voevodsky, Beilinson, Weil, Tate, Parshin, Kimura, Schur, and others, were recently extended from the realm of algebraic geometry to the broad noncommutative setting of differential graded (=dg) categories. This noncommutative viewpoint led to a proof of these celebrated conjectures in several new cases. Moreover, it enabled a proof of the noncommutative counterparts of the celebrated conjectures in many interesting cases. The purpose of this survey, written for a broad mathematical audience, is to give a rigorous overview of these recent developments.

\smallbreak

{\it Notations.}
Given a perfect base field $k$ of characteristic $p>0$, we will write $W(k)$ for its ring of $p$-typical Witt vectors and $K:=W(k)_{1/p}$ for the fraction field of $W(k)$. For example, when $k=\bbF_p$, we have $W(k)=\bbZ_p$ and $K=\bbQ_p$.
%-------------------------------------------------------------------------------
\section{Celebrated conjectures}\label{sec:conjectures}
%-------------------------------------------------------------------------------
In this section, we briefly recall some celebrated conjectures of Grothendieck, Voevodsky, Beilinson, Weil, Tate, Parshin, and Kimura (concerning smooth proper schemes), as well as a conjecture of Schur (concerning smooth schemes).
%-------------------------------------------------------------------------------
\subsection{Grothendieck standard conjecture of type $\mathrm{C}^+$}\label{sub:sign}
%-------------------------------------------------------------------------------
Let $k$ be a perfect base field of characteristic $p\geq 0$ and $X$ a smooth proper $k$-scheme of dimension $d$.  When $p=0$, we will write $H^\ast_{\mathrm{dR}}(X)$ for the de Rham cohomology of $X$. In the same vein, when $p>0$, we will write $H^\ast_{\mathrm{crys}}(X):=H^\ast_{\mathrm{crys}}(X/W(k))\otimes_{W(k)}K$ for the crystalline cohomology of $X$. Given an integer $0 \leq i \leq 2d$, consider the associated $i^{\mathrm{th}}$ K\"unneth projector $\pi^i\colon H^\ast_{\mathrm{dR}}(X) \to H^\ast_{\mathrm{dR}}(X)$, resp. $\pi^i\colon H^\ast_{\mathrm{crys}}(X) \to H^\ast_{\mathrm{crys}}(X)$, in de Rham cohomology, resp. in crystalline cohomology. In the sixties, Grothendieck \cite{Grothendieck} (see also \cite{Kleiman, Kleiman1}) conjectured the following:

\smallskip

{\it Conjecture $\mathrm{C}^+(X)$: The even K\"unneth projector $\pi^+:=\sum_{i\,\mathrm{even}} \pi^{i}$ is algebraic\footnote{If $\pi^+$  is algebraic, then the odd K\"unneth projector $\pi^-:=\sum_{i\,\mathrm{odd}} \pi^{i}$ is also algebraic.}.}

\smallskip

This conjecture is also usually called the ``sign conjecture''. It holds when $d\leq 2$, when $X$ is an abelian variety (see Kleiman \cite{Kleiman1}), and also when $k$ is a finite field (see Katz-Messing \cite{KM}). Besides these cases (and some other cases scattered in the literature), it remains wide open. 
\begin{remark}\label{rk:product}
Given smooth proper $k$-schemes $X$ and $Y$, we have the implication of conjectures $\mathrm{C}^+(X) + \mathrm{C}^+(Y) \Rightarrow \mathrm{C}^+(X\times Y)$.
\end{remark}
%-------------------------------------------------------------------------------
\subsection{Grothendieck standard conjecture of type $\mathrm{D}$}
%-------------------------------------------------------------------------------
Let $k$ be a perfect base field of characteristic $p\geq 0$ and $X$ a smooth proper $k$-scheme of dimension $d$. Consider the graded $\bbQ$-vector space $\cZ^\ast(X)_\bbQ/\!_{\sim\mathrm{hom}}$ of algebraic cycles on $X$ up to homological equivalence (when $p=0$, resp. $p>0$, we make use of de Rham cohomology, resp. crystalline cohomology). Consider also the graded $\bbQ$-vector space $\cZ^\ast(X)_\bbQ/\!_{\sim\mathrm{num}}$ of algebraic cycles on $X$ up to numerical equivalence. In the sixties, Grothendieck \cite{Grothendieck} (see also \cite{Kleiman, Kleiman1}) conjectured the following:

\smallskip

{\it Conjecture $\mathrm{D}(X)$: The equality $\cZ^\ast(X)_\bbQ/\!_{\sim\mathrm{hom}}=\cZ^\ast(X)_\bbQ/\!_{\sim\mathrm{num}}$ holds.}

\smallskip

This conjecture holds when $d\leq 2$, when $d\leq 4$ and $p=0$ (see Lieberman \cite{Lieberman}), and also when $X$ is an abelian variety and $p=0$ (see Lieberman \cite{Lieberman}). Besides these cases (and some other cases scattered in the literature), it remains wide open.

%-------------------------------------------------------------------------------
\subsection{Voevodsky nilpotence conjecture}
%-------------------------------------------------------------------------------
Let $k$ be a base field of characteristic $p\geq 0$ and $X$ a smooth proper $k$-scheme of dimension $d$. Following Voevodsky \cite{Voevodsky}, consider the graded $\bbQ$-vector space $\cZ^\ast(X)_\bbQ/\!_{\sim\mathrm{nil}}$ of algebraic cycles on $X$ up to nilpotence equivalence. In the nineties, Voevodsky \cite{Voevodsky} conjectured the following:

\smallskip

{\it Conjecture $\mathrm{V}(X)$: The equality $\cZ^\ast(X)_\bbQ/\!_{\sim\mathrm{nil}}=\cZ^\ast(X)_\bbQ/\!_{\sim\mathrm{num}}$ holds.}

\smallskip

This conjecture holds when $d\leq 2$ (see Voevodsky \cite{Voevodsky} and Voisin \cite{Voisin}), and also when $X$ is an abelian threefold and $p=0$ (see Kahn-Sebastian \cite{KS}). Besides these cases (and some other cases scattered in the literature), it remains wide open.
\begin{remark}\label{rk:implication}
Every algebraic cycle which is nilpotently trivial is also homologically trivial. Hence, we have the implication of conjectures $\mathrm{V}(X) \Rightarrow \mathrm{D}(X)$.
\end{remark}
%-------------------------------------------------------------------------------
\subsection{Beilinson conjecture}\label{sub:Beilinson}
%-------------------------------------------------------------------------------
Let $k=\bbF_q$ be a finite field of characteristic $p$ and $X$ a smooth proper $k$-scheme of dimension $d$. Consider the graded $\bbQ$-vector space $\cZ^\ast(X)_\bbQ/\!_{\sim\mathrm{rat}}$ of algebraic cycles on $X$ up to rational equivalence. In the eighties, Beilinson \cite{Beilinson} conjectured the following:

\smallskip

{\it Conjecture $\mathrm{B}(X)$: The equality $\cZ^\ast(X)_\bbQ/\!_{\sim\mathrm{rat}}=\cZ^\ast(X)_\bbQ/\!_{\sim\mathrm{num}}$ holds.}

\smallskip

This conjecture holds when $d\leq 1$, and also when $X$ is an abelian variety and $d\leq3$ (see Kahn \cite{Kahn1}). Besides these cases (and some other cases scattered in the literature), it remains wide open.

\begin{remark}\label{rk:Beilinson}
Every algebraic cycle which is rationally trivial is also nilpotently trivial. Hence, in the case where $k$ is a finite field, we have $\mathrm{B}(X) \Rightarrow \mathrm{V}(X)$.
\end{remark}

%-------------------------------------------------------------------------------
\subsection{Weil conjecture}\label{sub:Weil}
%-------------------------------------------------------------------------------
Let $k=\bbF_q$ be a finite field of characteristic $p$ and $X$ a smooth proper $k$-scheme of dimension $d$. Recall that the zeta function of $X$ is defined as the formal power series $Z(X;t) := \mathrm{exp}(\sum_{n \geq 1}\# X(\bbF_{q^n})\frac{t^n}{n}) \in \bbQ[\![t]\!]$, where $\mathrm{exp}(t):=\sum_{n \geq 0} \frac{t^n}{n!}$. In the same vein, given an integer $0 \leq i \leq 2d$, consider the formal power series $Z_i(X;t):=\mathrm{det}(\id - t \,\mathrm{Fr}^i |H^i_{\mathrm{crys}}(X))^{-1} \in K[\![t]\!]$, where $\mathrm{Fr}$ stands for the Frobenius endomorphism of $X$ and $\mathrm{Fr}^i$ for the induced automorphism of $H^i_{\mathrm{crys}}(X)$. Thanks to the Lefschetz trace formula established by Grothendieck and Berthelot (consult \cite{Berthelot}), we have the following weight decomposition:
\begin{equation}\label{eq:factorization}
Z(X;t)=\frac{Z_0(X;t) Z_2(X;t)\cdots Z_{2d}(X;t)}{Z_1(X;t)Z_3(X;t) \cdots Z_{2d-1}(X;t)} \in K[\![t]\!]\,.
\end{equation}

In the late forties, Weil \cite{Weil} conjectured the following\footnote{The above conjecture $\mathrm{W}(X)$ is a modern formulation of Weil's original conjecture; in the late forties crystalline cohomology was not yet developed.}:

\smallskip

{\it Conjecture $\mathrm{W}(X)$: The eigenvalues of the automorphism $\mathrm{Fr}^i$, with $0 \leq i \leq 2d$, are algebraic numbers and all their complex conjugates have absolute value $q^{\frac{i}{2}}$.}

\smallskip

In the particular case of curves, this famous conjecture follows from Weil's pioneering work \cite{Weil1}. Later, in the seventies, it was proved in full generality~by~Deligne\footnote{Deligne worked with \'etale cohomology instead. However, as explained by Katz-Messing in \cite{KM}, Deligne's results hold similarly in crystalline cohomology. More recently, Kedlaya \cite{Kedlaya} gave an alternative proof of the Weil conjecture which uses solely $p$-adic techniques.} \cite{Deligne}. In contrast with Weil's proof, which uses solely the classical intersection theory of divisors on surfaces, Deligne's proof makes use of several involved tools such as the theory of monodromy of Lefschetz pencils. The Weil conjecture has numerous applications. For example, when combined with the weight decomposition \eqref{eq:factorization}, it implies that the polynomials $\mathrm{det}(\id - t \,\mathrm{Fr}^i |H^i_{\mathrm{crys}}(X))$ have integer coefficients.

Recall that the Hasse-Weil zeta function of $X$ is defined as the (convergent) infinite product $\zeta(X;s):= \prod_{x \in X^0} (1- (q^{\mathrm{deg}(x)})^{-s})^{-1}$, with $\mathrm{Re}(s)>d$, where $X^0$ stands for the set of closed points of $X$ and $\mathrm{deg}(x)$ for the degree of the finite field extension $\kappa(x)/\bbF_q$. In the same vein, given an integer $0 \leq i \leq 2d$, consider the function $\zeta_i(X;s):= \mathrm{det}(\id - q^{-s} \,\mathrm{Fr}^i |H^i_{\mathrm{crys}}(X))^{-1}$. It follows from the Weil conjecture that $\zeta(X;s)=Z(X;q^{-s})$, with $\mathrm{Re}(s)>d$, and that $\zeta_i(X;s)=Z_i(X;q^{-s})$, with $\mathrm{Re}(s)>\frac{i}{2}$. Thanks to \eqref{eq:factorization}, we hence obtain the weight decomposition:
\begin{eqnarray}\label{eq:factorization1}
\zeta(X;s)= \frac{\zeta_0(X;s)\zeta_2(X;s) \cdots \zeta_{2d}(X;s)}{\zeta_1(X;s) \zeta_3(X;s) \cdots \zeta_{2d-1}(X;s)} && \mathrm{Re}(s)>d\,.
\end{eqnarray}
Note that \eqref{eq:factorization1} implies automatically that the Hasse-Weil zeta function $\zeta(X;s)$ of $X$ admits a (unique) meromorphic continuation to the entire complex plane.
\begin{remark}[Analogue of the Riemann hypothesis]\label{rk:RH}
The above conjecture $\mathrm{W}(X)$ is usually  called the ``analogue of the Riemann hypothesis'' because it implies that if $z\in \bbC$ is a pole of $\zeta_i(X;s)$, then $\mathrm{Re}(z)=\frac{i}{2}$. Consequently, if $z \in \bbC$ is a pole, resp. zero, of $\zeta(X;s)$, then $\mathrm{Re}(z)\in\{0, 1, \ldots, d\}$, resp. $\mathrm{Re}(z)\in\{\frac{1}{2}, \frac{2}{3}, \ldots, \frac{2d-1}{2}\}$. 
\end{remark}

%-------------------------------------------------------------------------------
\subsection{Tate conjecture}\label{sub:Tate}
%-------------------------------------------------------------------------------
Let $k=\bbF_q$ be a finite field of characteristic $p$ and $X$ a smooth proper $k$-scheme of dimension $d$. Given a prime number $l\neq p$, consider the associated $\bbQ_l$-linear cycle class map with values in $l$-adic cohomology:
\begin{equation}\label{eq:cycle}
\cZ^\ast(X)_{\bbQ_l}/\!_{\sim\mathrm{rat}} \too H^{2\ast}_{l\text{-}\mathrm{adic}}(X_{\overline{k}}, \bbQ_l(\ast))^{\mathrm{Gal}(\overline{k}/k)}\,.
\end{equation}
In the sixties, Tate \cite{Tate} conjectured the following:

\smallskip

{\it Conjecture $\mathrm{T}^l(X)$: The cycle class map \eqref{eq:cycle} is surjective.}

\smallskip

This conjecture holds when $d \leq 1$, when $X$ is an abelian variety and $d\leq 3$, and also when $X$ is a $K3$-surface; consult Totaro's survey \cite{Totaro}. Besides these cases (and some other cases scattered in the literature), it remains wide open.
%-------------------------------------------------------------------------------
\subsection{$p$-version of the Tate conjecture}
%-------------------------------------------------------------------------------
Let $k=\bbF_q$ be a finite field of characteristic $p$ and $X$ a smooth proper $k$-scheme of dimension $d$. Consider the associated $K$-linear cycle class map with values in crystalline cohomology (see \S\ref{sub:Weil}):
\begin{equation}\label{eq:cycle1}
\cZ^\ast(X)_K/_{\!\sim\mathrm{rat}} \too H^{2\ast}_{\mathrm{crys}}(X)(\ast)^{\mathrm{Fr}^{2\ast}}\,.
\end{equation}
Following Milne \cite{Milne}, the Tate conjecture admits the following $p$-version:

\smallskip

{\it Conjecture $\mathrm{T}^p(X)$: The cycle class map \eqref{eq:cycle1} is surjective.}

\smallskip

This conjecture is equivalent to $\mathrm{T}^l(X)$ (for every $l\neq p$) when $d\leq 3$. Hence, it also holds in the cases mentioned in \S\ref{sub:Tate}. Besides these cases (and some other cases scattered in the literature), it remains wide open. 
\begin{remark}
The $p$-version of the Tate conjecture can be alternatively formulated as follows: the $\bbQ_p$-linear cycle class map $\cZ^\ast(X)_{\bbQ_p}/_{\!\sim\mathrm{rat}} \to H^{2\ast}_{\mathrm{crys}}(X)(\ast)^{\mathrm{Fr}_p^{2\ast}}$, where $\mathrm{Fr}^{2\ast}_p$ stands for the crystalline Frobenius, is surjective.
\end{remark}
%-------------------------------------------------------------------------------
\subsection{Strong form of the Tate conjecture}
%-------------------------------------------------------------------------------
Let $k=\bbF_q$ be a finite field of characteristic $p$ and $X$ a smooth proper $k$-scheme of dimension $d$. Recall from \S\ref{sub:Weil} that the Hasse-Weil zeta function of $X$ is defined as the (convergent) infinite product $\zeta(X;s):=\prod_{x \in X^0}(1-(q^{\mathrm{deg}(x)})^{-s})^{-1}$, with $\mathrm{Re}(s)>d$. Moreover, as mentioned in {\em loc. cit.}, $\zeta(X;s)$ admits a meromorphic continuation to the entire complex plane. In the sixties, Tate \cite{Tate} also conjectured the following:

\smallskip

{\it Conjecture $\mathrm{ST}(X)$: The order $\mathrm{ord}_{s=j}\zeta(X;s)$ of the Hasse-Weil zeta function $\zeta(X;s)$ at the pole $s=j$, with $0\leq j \leq d$, is equal to $-\mathrm{dim}_\bbQ\cZ^j(X)_\bbQ/\!_{\sim\mathrm{num}}$.}

\begin{remark}\label{rk:strongTate}
As proved by Tate in \cite{Tate1}, resp. by Milne in \cite{Milne}, we have the equivalence of conjectures $\mathrm{ST}(X) \Leftrightarrow \mathrm{B}(X) + \mathrm{T}^l(X)$, resp. $\mathrm{ST}(X) \Leftrightarrow \mathrm{B}(X) + \mathrm{T}^p(X)$.
\end{remark}
Thanks to Remark \ref{rk:strongTate}, the conjecture $\mathrm{ST}(X)$ holds when $d \leq 1$, and also when $X$ is an abelian variety and $d\leq 3$. Besides these cases (and some other cases scattered in the literature), it remains wide open.
%-------------------------------------------------------------------------------
\subsection{Parshin conjecture}
%-------------------------------------------------------------------------------
Let $k=\bbF_q$ be a finite field of characteristic $p$ and $X$ a smooth proper $k$-scheme of dimension $d$. Consider the associated algebraic $K$-theory groups $K_n(X), n \geq 0$. In the eighties, Parshin conjectured the following:

\smallskip

{\it Conjecture $\mathrm{P}(X)$: The groups $K_n(X)$, with $n \geq 1$, are torsion.}

\smallskip

This conjecture holds when $d\leq 1$ (see Quillen \cite{Quillen} and Harder \cite{Harder}). Besides these cases (and some other cases scattered in the literature), it remains wide open.

\begin{remark}
As proved by Geisser in \cite{Geisser}, we have the implication of conjectures $\mathrm{B}(X) + \mathrm{ST}(X) \Rightarrow \mathrm{P}(X)$. This implies, in particular, that the conjecture $\mathrm{P}(X)$ also holds when $X$ is an abelian variety and $d \leq 3$.
\end{remark}
%-------------------------------------------------------------------------------
\subsection{Kimura-finiteness conjecture}
%-------------------------------------------------------------------------------
Let $k$ be a base field of characteristic $p\geq 0$ and $X$ a smooth proper $k$-scheme of dimension $d$. Consider the category of Chow motives $\Chow(k)_\bbQ$ introduced by Manin in \cite{Manin}. By construction, this category is $\bbQ$-linear, idempotent complete, symmetric monoidal, and comes equipped with a symmetric monoidal functor $\mathfrak{h}(-)_\bbQ\colon \mathrm{SmProp}(k)^\op \to \Chow(k)_\bbQ$ defined on smooth proper $k$-schemes. A decade ago, Kimura \cite{Kimura} conjectured the following:

\smallskip

{\it Conjecture $\mathrm{K}(X)$: The Chow motive $\mathfrak{h}(X)_\bbQ$ is Kimura-finite\footnote{Let $(\cC,\otimes, {\bf 1})$ be a $\bbQ$-linear, idempotent complete, symmetric monoidal category. Following Kimura \cite{Kimura}, recall that an object $a \in \cC$ is called {\em even-dimensional}, resp. {\em odd-dimensional}, if $\wedge^n(a)\simeq 0$, resp. $\mathrm{Sym}^n(a)\simeq 0$, for some $n\gg 0$. The biggest integer $\mathrm{kim}_+(a)$, resp. $\mathrm{kim}_-(a)$, for which $\wedge^{\mathrm{kim}_+(a)}(a) \not\simeq 0$, resp. $\mathrm{Sym}^{\mathrm{kim}_-(a)}(a)\not\simeq 0$, is called the {\em even Kimura-dimension}, resp. {\em odd Kimura-dimension}, of $a$. Recall also that an object $a\in \cC$ is called {\em Kimura-finite} if $a\simeq a_+\oplus a_-$, with $a_+$ even-dimensional and $a_-$ odd-dimensional. The integer $\mathrm{kim}(a):=\mathrm{kim}_+(a_+) + \mathrm{kim}_-(a_-)$ is called the {\em Kimura-dimension} of $a$.}.}

\smallskip
 
This conjecture holds when $d\leq 1$ and also when $X$ is an abelian variety (see Kimura \cite{Kimura} and Shermenev \cite{Shermenev}). Besides these cases\footnote{In the particular case where $k=\overline{k}$, $p=0$, and $X$ is a surface with $p_g(X)=0$, Guletskii and Pedrini proved in \cite{GP} that the conjecture $\mathrm{K}(X)$ is equivalent to a celebrated conjecture of Bloch \cite{Bloch} concerning the vanishing of the Albanese kernel.} (and some other cases scattered in the literature), it remains wide open.
%-------------------------------------------------------------------------------
\subsection{Schur-finiteness conjecture}
%-------------------------------------------------------------------------------
Let $k$ be a perfect base field of characteristic $p\geq 0$ and $X$ a smooth $k$-scheme of dimension $d$. Consider the triangulated category of geometric mixed motives $\mathrm{DM}_{\mathrm{gm}}(k)_\bbQ$ introduced by Voevodsky in \cite{Voevodsky2}. By construction, this category is $\bbQ$-linear, idempotent complete, symmetric monoidal, and comes equipped with a symmetric monoidal functor $M(-)_\bbQ\colon \mathrm{Sm}(k) \to \mathrm{DM}_{\mathrm{gm}}(k)_\bbQ$ defined on smooth $k$-schemes. Moreover, as proved in {\em loc. cit.}, the classical category of Chow motives $\Chow(k)_\bbQ$ may be embedded fully-faithfully into $\mathrm{DM}_{\mathrm{gm}}(k)_\bbQ$. An important conjecture in the theory of mixed motives is the following:

\smallskip

{\it Conjecture $\mathrm{S}(X)$: The mixed motive $M(X)_\bbQ$ is Schur-finite\footnote{Let $(\cC,\otimes, {\bf 1})$ be a $\bbQ$-linear, idempotent complete, symmetric monoidal category. Following Deligne \cite{Deligne-Moscow}, every partition $\lambda$ of an integer $n\geq 1$ gives naturally rise to a Schur-functor $S_\lambda\colon \cC \to \cC$. For example, when $\lambda=(1, \ldots, 1)$, resp. $\lambda=(n)$, we have $S_{(1, \ldots, 1)}(a)= \wedge^n(a)$, resp. $S_{(1)}(a)=\mathrm{Sym}^n(a)$. An object $a \in \cC$ is called {\em Schur-finite} if $S_\lambda(a)\simeq 0$~for~some~partition~$\lambda$.}.}

\smallskip

This conjecture holds when $d\leq 1$ (see Guletskii \cite{Guletskii} and Mazza \cite{Mazza}) and also when $X$ is an abelian variety (see Kimura \cite{Kimura} and Shermenev \cite{Shermenev}). Besides these cases (and some other cases scattered in the literature), it remains wide open.

\begin{remark}
It is well-known that Kimura-finiteness implies Schur-finiteness. However, the converse does {\em not} holds. For example, O'Sullivan constructed a certain smooth surface $X$ whose mixed motive $M(X)_\bbQ$ is Schur-finite but {\em not} Kimura-finite; consult \cite{Mazza}. An important open problem is the classification of all the Kimura-finite mixed motives and the computation of the corresponding Kimura-dimensions.
\end{remark}
%-------------------------------------------------------------------------------
\section{Noncommutative counterparts}\label{sec:NCcounterparts}
%-------------------------------------------------------------------------------
In this section we describe the noncommutative counterparts of the celebrated conjectures of \S\ref{sec:conjectures}.
We will assume some basic familiarity with the language of differential graded (=dg) categories; consult Keller's survey \cite{Keller}. In particular, we will use freely the notion of {\em smooth proper} dg category in the sense of Kontsevich \cite{Miami, finMot, Lefschetz, Kontsevich-talk, IAS}. Examples include the finite-dimensional algebras of finite global dimension $A$ (over a perfect base field) as well as the dg categories of perfect complexes $\perf_\dg(X)$ associated to smooth proper schemes $X$ (or, more generally, to smooth proper algebraic stacks $\cX$). In addition, we will make essential use of the recent theory of noncommutative motives; consult the book \cite{book} and~the~survey~\cite{survey}. 
%-------------------------------------------------------------------------------
\subsection{Noncommutative Grothendieck standard conjecture of type $\mathrm{C}^+$}\label{sub:standardC}
%-------------------------------------------------------------------------------
Let $k$ be a perfect base field of characteristic $p\geq 0$ and $\cA$ a smooth proper $k$-linear dg category. In what follows, we will write $\dgcat_{\mathrm{sp}}(k)$ for the category of (essentially small) $k$-linear dg categories. Recall from \cite[\S9]{JEMS} that, when $p=0$, periodic cyclic homology gives rise to a symmetric monoidal functor  
\begin{equation}\label{eq:HP}
HP_\pm(-)\colon \dgcat_{\mathrm{sp}}(k) \too \mathrm{vect}_{\bbZ/2}(k)
\end{equation}
with values in the category of finite-dimensional $\bbZ/2$-graded $k$-vector spaces. In the same vein, recall from \cite[\S2]{NCpositive} that, when $p>0$, topological periodic cyclic homology\footnote{Topological periodic cyclic homology is defined as the Tate cohomology of the circle group action on topological Hochschild homology; consult Hesselholt \cite{Lars} and Nikolaus-Scholze \cite{Scholze}.} gives rise to a symmetric monoidal functor
\begin{equation}\label{eq:TP}
TP_\pm(-)_{1/p}\colon \dgcat_{\mathrm{sp}}(k) \too \mathrm{vect}_{\bbZ/2}(K)
\end{equation}
with values in the category of finite-dimensional $\bbZ/2$-graded $K$-vector spaces. 
\begin{remark}[Relation with de Rham cohomology and crystalline cohomology]
The above functor \eqref{eq:HP}, resp. \eqref{eq:TP}, may be understood as the noncommutative counterpart of de Rham cohomology, resp. crystalline cohomology. Concretely, given a smooth proper $k$-scheme $X$, we have the following natural isomorphisms of finite-dimensional $\bbZ/2$-graded vector spaces:
\begin{eqnarray}
HP_\pm(\perf_\dg(X))  & \simeq &  (\bigoplus_{i\,\mathrm{even}}H^i_{\mathrm{dR}}(X), \bigoplus_{i\,\mathrm{odd}}H^i_{\mathrm{dR}}(X)) \label{eq:iso1}\\
TP_\pm(\perf_\dg(X))_{1/p}  & \simeq & (\bigoplus_{i\,\mathrm{even}}H^i_{\mathrm{crys}}(X), \bigoplus_{i\,\mathrm{odd}}H^i_{\mathrm{crys}}(X))\,. \label{eq:iso2}
\end{eqnarray}
On the one hand, \eqref{eq:iso1} follows from the classical Hochschild-Kostant-Rosenberg theorem; see Feigin-Tsygan \cite{FT}. On the other hand, \eqref{eq:iso2} follows from the recent work of Scholze on integral $p$-adic Hodge theory; consult \cite{Elment}\cite[Thm.~5.2]{NCpositive-CD}.
\end{remark}
Recall from \cite[\S4.1]{book} the definition of the category of noncommutative Chow motives $\NChow(k)_\bbQ$. By construction, this category is $\bbQ$-linear, idempotent complete, symmetric monoidal, and comes equipped with a symmetric monoidal functor $U(-)_\bbQ\colon \dgcat_{\mathrm{sp}}(k) \to \NChow(k)_\bbQ$. Moreover, we have a natural isomorphism
\begin{equation}\label{eq:Homs}
\Hom_{\NChow(k)_\bbQ}(U(k)_\bbQ, U(\cA)_\bbQ)\simeq K_0(\cD_c(\cA))_\bbQ=:K_0(\cA)_\bbQ\,,
\end{equation}
where $\cD(\cA)$ stands for the derived category of $\cA$ and $\cD_c(\cA)$ for its full triangulated subcategory of compact objects. As proved in \cite[Thm.~9.2]{JEMS} when $p=0$, resp. in \cite[Thm.~2.3]{NCpositive} when $p>0$, the above functor \eqref{eq:HP}, resp. \eqref{eq:TP}, descends to the category of noncommutative Chow motives. 

Consider the even K\"unneth projector 
\begin{eqnarray*}
\pi_+\colon HP_\pm(\cA) \to HP_\pm(\cA) && \text{resp.}\,\,\,\pi_+\colon TP_\pm(\cA)_{1/p} \to TP_\pm(\cA)_{1/p}
\end{eqnarray*}
in periodic cyclic homology, resp. in topological periodic cyclic homology. This projector is {\em algebraic} if there exists an endomorphism $\underline{\pi}_+\colon U(\cA)_\bbQ \to U(\cA)_\bbQ$ such that $HP_\pm(\underline{\pi}_+)=\pi_+$, resp. $TP_\pm(\underline{\pi}_+)_{1/p}=\pi_+$. Under these definitions,~the~Grothendieck standard conjecture of type $\mathrm{C}^+$ admits the following noncommutative counterpart: 

\smallskip

{\it Conjecture $\mathrm{C}_{\mathrm{nc}}^+(\cA)$: The even K\"unneth projector $\pi_+$ is algebraic\footnote{If $\pi_+$  is algebraic, then the odd K\"unneth projector $\pi_-$ is also algebraic.}.}

\begin{remark}
Similarly to Remark \ref{rk:product}, given smooth proper $k$-linear dg categories $\cA$ and $\cB$, we have the implication of conjectures $\mathrm{C}_{\mathrm{nc}}^+(\cA) + \mathrm{C}_{\mathrm{nc}}^+(\cB) \Rightarrow \mathrm{C}_{\mathrm{nc}}^+(\cA\otimes \cB)$.
\end{remark}

The next result relates this conjecture with Grothendieck's original conjecture:

\begin{theorem}{(\cite[Thm.~1.1]{NC-zero-CD} and \cite[Thm.~1.1]{NCpositive-CD})}\label{thm:main1}
Given a smooth proper $k$-scheme $X$, we have the equivalence of conjectures $\mathrm{C}^+(X) \Leftrightarrow \mathrm{C}^+_{\mathrm{nc}}(\perf_\dg(X))$.
\end{theorem} 
%-------------------------------------------------------------------------------
\subsection{Noncommutative Grothendieck standard conjecture of type $\mathrm{D}$}
%-------------------------------------------------------------------------------
Let $k$ be a perfect base field of characteristic $p\geq0$ and $\cA$ a smooth proper $k$-linear dg category. Note that by combining the above isomorphism \eqref{eq:Homs} with the functor \eqref{eq:HP}, resp. \eqref{eq:TP}, we obtain an induced $\bbQ$-linear homomorphism:
\begin{eqnarray}\label{eq:induced}
K_0(\cA)_\bbQ \too HP_+(\cA) && \text{resp.}\,\,\,K_0(\cA)_\bbQ \too TP_+(\cA)_{1/p}\,.
\end{eqnarray}
The homomorphism \eqref{eq:induced} may be understood as the noncommutative counterpart of the cycle class map. In what follows, we will write $K_0(\cA)_\bbQ/\!_{\sim\mathrm{hom}}$ for the quotient of $K_0(\cA)_\bbQ$ by the kernel of $\eqref{eq:induced}$. Consider also the Euler bilinear pairing:
\begin{eqnarray*}
\chi\colon K_0(\cA)\times K_0(\cA) \too \bbZ && ([M],[N])\mapsto \sum_{n\in \bbZ} (-1)^n \mathrm{dim}_k\Hom_{\cD_c(\cA)}(M,N[n])\,.
\end{eqnarray*}
This pairing is not symmetric neither skew-symmetric. Nevertheless, as proved in \cite[Prop.~4.24]{book}, the left and right kernels of $\chi$ agree. In what follows, we will write $K_0(\cA)/\!_{\sim\mathrm{num}}$ for the quotient of $K_0(\cA)$ by the kernel of $\chi$ and $K_0(\cA)_\bbQ/\!_{\sim\mathrm{num}}$ for the $\bbQ$-vector space\footnote{As proved in \cite[Thm.~5.1]{NCpositive}, $K_0(\cA)/\!_{\sim\mathrm{num}}$ is a finitely generated free abelian group. Consequently, $K_0(\cA)_\bbQ/\!_{\sim\mathrm{num}}$ is a finite-dimensional $\bbQ$-vector space.} $K_0(\cA)/\!_{\sim\mathrm{num}}\otimes_\bbZ \bbQ$. Under these definitions, the Grothendieck standard conjecture of type $\mathrm{D}$ admits the following noncommutative counterpart:

\smallskip

{\it Conjecture $\mathrm{D}_{\mathrm{nc}}(\cA)$: The equality $K_0(\cA)_\bbQ/\!_{\sim\mathrm{hom}} = K_0(\cA)_\bbQ/\!_{\sim\mathrm{num}}$ holds.}

\smallskip

The next result relates this conjecture with Grothendieck's original conjecture:

\begin{theorem}{(\cite[Thm.~1.1]{NC-zero-CD} and \cite[Thm.~1.1]{NCpositive-CD})}\label{thm:main2}
Given a smooth proper $k$-scheme $X$, we have the equivalence of conjectures $\mathrm{D}(X) \Leftrightarrow \mathrm{D}_{\mathrm{nc}}(\perf_\dg(X))$.
\end{theorem}
%-------------------------------------------------------------------------------
\subsection{Noncommutative Voevodsky nilpotence conjecture}
%-------------------------------------------------------------------------------
Let $k$ be a base field of characteristic $p\geq 0$ and $\cA$ a smooth proper $k$-linear dg category. Similarly to Voevodsky's definition of the nilpotence equivalence relation, an element $\alpha$ of the Grothendieck group $K_0(\cA)_\bbQ$ is called {\em nilpotently trivial} if there exists an integer $n \gg0$ such that the associated element $\alpha^{\otimes n}$ of the Grothendieck group $K_0(\cA^{\otimes n})_\bbQ$ is equal to zero. In what follows, we will write $K_0(\cA)_\bbQ/\!_{\sim\mathrm{nil}}$ for the quotient of $K_0(\cA)_\bbQ$ by the nilpotently trivial elements. Under these definitions, the Voevodsky nilpotence conjecture admits the following noncommutative counterpart:

\smallskip

{\it Conjecture $\mathrm{V}_{\mathrm{nc}}(\cA)$: The equality $K_0(\cA)_\bbQ/\!_{\sim\mathrm{nil}} = K_0(\cA)_\bbQ/\!_{\sim\mathrm{num}}$ holds.}

\begin{remark}
The image of a nilpotently trivial element $\alpha \in K_0(\cA)_\bbQ$ under the above $\bbQ$-linear homomorphism \eqref{eq:induced} is equal to zero. Consequently, similarly to Remark \ref{rk:implication}, we have the implication of conjectures $\mathrm{V}_{\mathrm{nc}}(\cA) \Rightarrow \mathrm{D}_{\mathrm{nc}}(\cA)$. 
\end{remark}

The next result relates this conjecture with Voevodsky's original conjecture:

\begin{theorem}{(\cite[Thm.~1.1]{Crelle})}\label{thm:main3}
Given a smooth proper $k$-scheme $X$, we have the equivalence of conjectures $\mathrm{V}(X) \Leftrightarrow \mathrm{V}_{\mathrm{nc}}(\perf_\dg(X))$.
\end{theorem}
 %-------------------------------------------------------------------------------
\subsection{Noncommutative Beilinson conjecture}\label{sub:NCBeilinson}
%-------------------------------------------------------------------------------
Let $k=\bbF_q$ be a finite field of characteristic $p$ and $\cA$ a smooth proper $k$-linear dg category. The Beilinson conjecture admits the following noncommutative counterpart:

\smallskip

{\it Conjecture $\mathrm{B}_{\mathrm{nc}}(\cA)$: The equality $K_0(\cA)_\bbQ = K_0(\cA)_\bbQ/\!_{\sim\mathrm{num}}$ holds.}

\begin{remark}
Similarly to Remark \ref{rk:Beilinson}, note that in the case where $k$ is a finite field, we have the implication of conjectures $\mathrm{B}_{\mathrm{nc}}(\cA) \Rightarrow \mathrm{V}_{\mathrm{nc}}(\cA)$.
\end{remark}

The next result relates this conjecture with Beilinson's original conjecture:

\begin{theorem}{(\cite[Thm.~1.3]{NCTate})}\label{thm:main4}
Given a smooth proper $k$-scheme $X$, we have the equivalence of conjectures $\mathrm{B}(X) \Leftrightarrow \mathrm{B}_{\mathrm{nc}}(\perf_\dg(X))$.
\end{theorem}
%-------------------------------------------------------------------------------
\subsection{Noncommutative Weil conjecture}\label{sub:NCWeil}
%-------------------------------------------------------------------------------
Let $k=\bbF_q$ be a finite field of characteristic $p$ and $\cA$ a smooth proper $k$-linear dg category. As explained in \cite[\S6]{NCWeil}, the topological periodic cyclic homology group $TP_0(\cA)_{1/p}$, resp. $TP_1(\cA)_{1/p}$, comes equipped with an automorphism $\mathrm{F}_0$, resp. $\mathrm{F}_1$, called the ``cyclotomic Frobenius''\footnote{The cyclotomic Frobenius is not compatible with the $\bbZ/2$-graded structure of $TP_\ast(\cA)_{1/p}$. Instead, we have canonical isomorphisms $\mathrm{F}_n \simeq q\cdot \mathrm{F}_{n+2}$ for every $n \in \bbZ$.}. Hence, we define the {\em even/odd zeta function of $\cA$} as the formal power series:
\begin{eqnarray*}
Z_{\mathrm{even}}(\cA;t)& := & \mathrm{det}(\id - t\,\mathrm{F}_0|TP_0(\cA)_{1/p})^{-1} \in K[\![t]\!] \\
Z_{\mathrm{odd}}(\cA;t)& := & \mathrm{det}(\id - t\,\mathrm{F}_1|TP_1(\cA)_{1/p})^{-1} \in K[\![t]\!] \,.
\end{eqnarray*}
Under these definitions, Weil's conjecture admits the noncommutative counterpart:

\smallskip

{\it Conjecture $\mathrm{W}_{\mathrm{nc}}(\cA)$: The eigenvalues of the automorphism $\mathrm{F}_0$, resp. $\mathrm{F}_1$, are algebraic numbers and all their complex conjugates have absolute~value~$1$,~resp.~$q^{\frac{1}{2}}$.}

\smallskip

In contrast with the commutative world, the cyclotomic Frobenius is not induced from an endomorphism\footnote{Note that in the particular case where $\cA$ is a $k$-algebra $A$ the Frobenius map $a \mapsto a^q$ is a $k$-algebra endomorphism if and only if $A$ is commutative.} of $\cA$. Consequently, in contrast with the commutative world, it is not known if the polynomials $\mathrm{det}(\id - t\,\mathrm{F}_0|TP_0(\cA)_{1/p})$ and $\mathrm{det}(\id - t\,\mathrm{F}_1|TP_1(\cA)_{1/p})$ have integer coefficients (or rational coefficients). Nevertheless, after choosing an embedding $\iota\colon K\hookrightarrow \bbC$, we define the {\em even/odd Hasse-Weil zeta function of $\cA$} as follows:
\begin{eqnarray*}
\zeta_{\mathrm{even}}(\cA;s)& := & \mathrm{det}(\id - q^{-s} (\mathrm{F}_0 \otimes_{K, \iota}\bbC)\,|\,TP_0(\cA)_{1/p}\otimes_{K, \iota} \bbC)^{-1} \\
\zeta_{\mathrm{odd}}(\cA;s)& := & \mathrm{det}(\id - q^{-s} (\mathrm{F}_1 \otimes_{K, \iota}\bbC)\,|\,TP_1(\cA)_{1/p}\otimes_{K, \iota} \bbC)^{-1}\,.
\end{eqnarray*}
\begin{remark}[Analogue of the noncommutative Riemann hypothesis]\label{rk:key}
Similarly to Remark \ref{rk:RH}, the conjecture $\mathrm{W}_{\mathrm{nc}}(\cA)$ may be called the ``analogue of the noncommutative Riemann hypothesis'' because it implies that if $z \in \bbC$ is a pole of $\zeta_{\mathrm{even}}(\cA; s)$, resp. $\zeta_{\mathrm{odd}}(\cA; s)$, then $\mathrm{Re}(z)=0$, resp. $\mathrm{Re}(z)=\frac{1}{2}$ (independently of the chosen $\iota$).
\end{remark}
The next result relates the above conjecture with Weil's original conjecture:

\begin{theorem}{(\cite[Thm.~1.5]{NCWeil})}
Given a smooth proper $k$-scheme $X$, we have the equivalence of conjectures $\mathrm{W}(X) \Leftrightarrow \mathrm{W}_{\mathrm{nc}}(\perf_\dg(X))$.
\end{theorem}
%-------------------------------------------------------------------------------
\subsection{Noncommutative Tate conjecture}\label{sub:NCTate}
%-------------------------------------------------------------------------------
Let $k=\bbF_q$ be a finite field of characteristic $p$ and $\cA$ a smooth proper $k$-linear dg category. Given a prime number $l\neq p$, consider the following abelian groups
\begin{eqnarray}\label{eq:groups-Tate}
\Hom(\bbZ(l^\infty), \pi_{-1}(L_{KU}(K(\cA\otimes_k\bbF_{q^n})))) && n \geq 1\,,
\end{eqnarray}
where $\bbZ(l^\infty)$ stands for the Pr\"ufer $l$-group\footnote{The functor $\Hom(\bbZ(l^\infty),-)$ agrees with the classical $l$-adic Tate module functor.}, $K(\cA\otimes_k \bbF_{q^n})$ for the algebraic $K$-theory spectrum of the dg category $\cA\otimes_k \bbF_{q^n}$, and $L_{KU}(-)$ for the Bousfield localization functor with respect to topological complex $K$-theory $KU$. Under these notations, the Tate conjecture admits the following noncommutative counterpart:

\smallskip

{\it Conjecture $\mathrm{T}^l_{\mathrm{nc}}(\cA)$: The abelian groups \eqref{eq:groups-Tate} are trivial.}

\begin{remark}
Note that the conjecture $\mathrm{T}^l_{\mathrm{nc}}(\cA)$ holds, for example, whenever the abelian groups $\pi_{-1}(L_{KU}(K(\cA\otimes_k \bbF_{q^n}))), n\geq 1$, are finitely generated. 
\end{remark}
The next result, obtained by leveraging the pioneering work of Thomason \cite{Thomason-Tate}, relates this conjecture with Tate's original conjecture:
\begin{theorem}{(\cite[Thm.~1.3]{NCTate})}\label{thm:main5}
Given a smooth proper $k$-scheme $X$, we have the equivalence of conjectures $\mathrm{T}^l(X) \Leftrightarrow \mathrm{T}^l_{\mathrm{nc}}(\perf_\dg(X))$.
\end{theorem}
%-------------------------------------------------------------------------------
\subsection{Noncommutative $p$-version of the Tate conjecture}\label{sub:NCTate1}
%-------------------------------------------------------------------------------
Let $k=\bbF_q$ be a finite field of characteristic $p$ and $\cA$ a smooth proper $k$-linear dg category. Recall from \S\ref{sub:NCWeil} that the $K$-vector space $TP_0(\cA)_{1/p}$ comes equipped with an automorphism $\mathrm{F}_0$ called the ``cyclotomic Frobenius''. Moreover, as explained in \cite[\S3]{NCTate}, the right-hand side of \eqref{eq:induced} gives rise to a $K$-linear homomorphism:
\begin{equation}\label{eq:induced1}
K_0(\cA)_K \too TP_0(\cA)_{1/p}^{\mathrm{F}_0}\,.
\end{equation}
Under these notations, the $p$-version of the Tate conjecture admits the following noncommutative counterpart:

\smallskip

{\it Conjecture $\mathrm{T}^p_{\mathrm{nc}}(\cA)$: The homomorphism \eqref{eq:induced1} is surjective.}

\smallskip

The next result relates this conjecture with the original conjecture:
\begin{theorem}{(\cite[Thm.~1.3]{NCTate})}\label{thm:main6}
Given a smooth proper $k$-scheme $X$, we have the equivalence of conjectures $\mathrm{T}^p(X) \Leftrightarrow \mathrm{T}^p_{\mathrm{nc}}(\perf_\dg(X))$.
\end{theorem}
%-------------------------------------------------------------------------------
\subsection{Noncommutative strong form of the Tate conjecture}\label{sub:NCTate2}
%-------------------------------------------------------------------------------
Let $k=\bbF_q$ be a finite field of characteristic $p$ and $\cA$ a smooth proper $k$-linear dg category. Recall from \S\ref{sub:NCWeil} the definition of the even Hasse-Weil zeta function $\zeta_{\mathrm{even}}(\cA;s)$ of $\cA$. Under these notations, the strong of the Tate conjecture admits the following noncommutative counterpart:

\smallskip

{\it Conjecture $\mathrm{ST}_{\mathrm{nc}}(\cA)$: The order $\mathrm{ord}_{s=0}\zeta_{\mathrm{even}}(\cA;s)$ of the even Hasse-Weil zeta function $\zeta_{\mathrm{even}}(\cA;s)$ at the pole $s=0$ is equal to $-\mathrm{dim}_\bbQ K_0(\cA)_\bbQ/\!_{\sim\mathrm{num}}$.}
\begin{remark}[Alternative formulation]\label{rk:multiplicity}
By definition of the even Hasse-Weil zeta function of $\cA$, the integer $-\mathrm{ord}_{s=0}\zeta_{\mathrm{even}}(\cA;s)$ agrees with the algebraic multiplicity of the eigenvalue $q^0=1$ of the automorphism $\mathrm{F}_0\otimes_{K, \iota}\bbC$ (or, equivalently, of $\mathrm{F}_0$). Hence, the conjecture $\mathrm{ST}_{\mathrm{nc}}(\cA)$ may be alternatively formulated as follows: {\em the algebraic multiplicity of the eigenvalue $1$ of $\mathrm{F}_0$ agrees with $\mathrm{dim}_\bbQ K_0(\cA)_\bbQ/_{\!\sim \mathrm{num}}$}. This shows, in particular, that the integer $\mathrm{ord}_{s=0}\zeta_{\mathrm{even}}(\cA;s)$ is independent of the embedding $\iota\colon K \hookrightarrow \bbC$ used in the definition of $\zeta_{\mathrm{even}}(\cA;s)$.
\end{remark}

\begin{remark}
Similarly to Remark \ref{rk:strongTate}, as proved in \cite[Thm.~9.3]{NCWeil}, we have the equivalence of conjectures $\mathrm{ST}_{\mathrm{nc}}(\cA) \Leftrightarrow \mathrm{B}_{\mathrm{nc}}(\cA) + \mathrm{T}^p_{\mathrm{nc}}(\cA)$.
\end{remark}

The next result relates this conjecture with Tate's original conjecture:
\begin{theorem}{(\cite[Thm.~1.17]{NCWeil})}\label{thm:main7}
Given a smooth proper $k$-scheme $X$, we have the equivalence of conjectures $\mathrm{ST}(X) \Leftrightarrow \mathrm{ST}_{\mathrm{nc}}(\perf_\dg(X))$.
\end{theorem}

%-------------------------------------------------------------------------------
\subsection{Noncommutative Parshin conjecture}
%-------------------------------------------------------------------------------
Let $k=\bbF_q$ be a finite field of characteristic $p$ and $\cA$ a smooth proper $k$-linear dg category. The Parshin conjecture admits the following noncommutative counterpart:

\smallskip

{\it Conjecture $\mathrm{P}_{\mathrm{nc}}(\cA)$: The groups $K_n(\cA)$, with $n \geq 1$, are torsion.}

\smallskip

The next result relates this conjecture with Parshin's original conjecture:
\begin{theorem}{(\cite[Thm.~1.3]{NCTate})}\label{thm:main8}
Given a smooth proper $k$-scheme $X$, we have the equivalence of conjectures $\mathrm{P}(X) \Leftrightarrow \mathrm{P}_{\mathrm{nc}}(\perf_\dg(X))$.
\end{theorem}
%-------------------------------------------------------------------------------
\subsection{Noncommutative Kimura-finiteness conjecture}
%-------------------------------------------------------------------------------
Let $k$ be a base field of characteristic $p\geq 0$ and $\cA$ a smooth proper $k$-linear dg category. Recall from \S\ref{sub:standardC} that, by construction, the category of noncommutative Chow motives $\NChow(k)_\bbQ$ is $\bbQ$-linear, idempotent complete and symmetric monoidal. Hence, the Kimura-finiteness conjecture admits the following noncommutative counterpart:

\smallskip

{\it Conjecture $\mathrm{K}_{\mathrm{nc}}(\cA)$: The noncommutative Chow motive $U(\cA)_\bbQ$ is Kimura-finite.}

\smallskip

The next result relates this conjecture with Kimura's original conjecture:
\begin{theorem}{(\cite[Thm.~2.1]{JNCG})}\label{thm:main9}
Given a smooth proper $k$-scheme $X$, we have the implication of conjectures $\mathrm{K}(X) \Rightarrow \mathrm{K}_{\mathrm{nc}}(\perf_\dg(X))$.
\end{theorem}

%-------------------------------------------------------------------------------
\subsection{Noncommutative Schur-finiteness conjecture}
%-------------------------------------------------------------------------------
Let $k$ be a perfect base field of characteristic $p\geq 0$ and $\cA$ a smooth $k$-linear dg category. Recall from \cite[\S8-\S9]{book} the definition of the triangulated category of noncommutative mixed motives $\mathrm{NMot}(k)_\bbQ$ (denoted by $\mathrm{Nmot}^{\bbA^1}_{\mathrm{loc}}(k)_\bbQ$ in {\em loc. cit.}). By construction, this category is $\bbQ$-linear, idempotent complete, symmetric monoidal, and comes equipped with a symmetric monoidal functor $\mathrm{U}(-)_\bbQ\colon \dgcat_{\mathrm{s}}(k)\to \mathrm{NMot}(k)_\bbQ$ defined on smooth dg categories. Under these notations, the Schur-finiteness conjecture admits the following noncommutative counterpart:

\smallskip

{\it Conjecture $\mathrm{S}_{\mathrm{nc}}(\cA)$: The noncommutative mixed motive $\mathrm{U}(\cA)_\bbQ$ is Schur-finite.} 

\smallskip

The next result relates this conjecture with Schur's original conjecture:
\begin{theorem}{(\cite[Prop.~9.17]{book})}\label{thm:main11}
Given a smooth $k$-scheme $X$, we have the equivalence of conjectures $\mathrm{S}(X) \Leftrightarrow \mathrm{S}_{\mathrm{nc}}(\perf_\dg(X))$.
\end{theorem}
%-------------------------------------------------------------------------------
\section{Applications to commutative geometry}\label{sec:applications1}
%-------------------------------------------------------------------------------
Morally speaking, the theorems of \S\ref{sec:NCcounterparts} show that the celebrated conjectures of Grothendieck, Voevodsky, Beilinson, Weil, Tate, Parshin, and Schur, belong not only to the realm of algebraic geometry but also to the broad setting of dg categories. This noncommutative viewpoint, where one studies a scheme via its dg category of perfect complexes, led to a proof\footnote{In what concerns the Weil conjecture (and the Grothendieck standard conjecture of type $\mathrm{C}^+$ over a finite field), the noncommutative viewpoint led to an alternative proof of this celebrated conjecture in several new cases, which avoids all the involved tools used by Deligne.} of these celebrated conjectures in several new cases. In this section, we describe some of these new cases. 
\begin{notation}
In order to simplify the exposition, we will often use the letter $\mathrm{C}$ to denote one of the celebrated conjectures $\{\mathrm{C}^+, \mathrm{D}, \mathrm{V}, \mathrm{B}, \mathrm{W}, \mathrm{T}^l, \mathrm{T}^p, \mathrm{ST}, \mathrm{P}, \mathrm{K}, \mathrm{S}\}$.
\end{notation}
%-------------------------------------------------------------------------------
\subsection{Derived invariance}
%-------------------------------------------------------------------------------
Let $k$ be a base field of characteristic $p\geq 0$. Note that the theorems of \S\ref{sec:NCcounterparts} imply automatically the following result:
\begin{corollary}\label{cor:derived}
Let $X$ and $Y$ be two smooth proper $k$-schemes (in the case of conjecture $\mathrm{S}$ we assume solely that $X$ and $Y$ are smooth) with (Fourier-Mukai) equivalent categories of perfect complexes $\perf(X)$ and $\perf(Y)$. Under these assumptions, we have the following equivalences of conjectures:
\begin{eqnarray*}
\mathrm{C}(X) \Leftrightarrow \mathrm{C}(Y) & \mathrm{with} & \mathrm{C} \in \{\mathrm{C}^+, \mathrm{D}, \mathrm{V}, \mathrm{B}, \mathrm{W}, \mathrm{T}^l, \mathrm{T}^p, \mathrm{ST}, \mathrm{P}, \mathrm{S}\}.
\end{eqnarray*}
\end{corollary}
Roughly speaking, Corollary \ref{cor:derived} shows that the celebrated conjectures of \S\ref{sec:conjectures} are invariant under derived equivalence. This flexibility is very useful and is often used in the proofs of some theorems below.
%-------------------------------------------------------------------------------
\subsection{Quadric fibrations}
%-------------------------------------------------------------------------------
Let $k$ be a perfect base field of characteristic $p\geq 0$, $B$ a smooth proper $k$-scheme of dimension $d$ (in the case of conjecture $\mathrm{S}$ we assume solely that $B$ is smooth), and $q\colon Q \to B$ a flat quadric fibration of relative dimension $d_q$.
\begin{theorem}\label{thm:app1-1}
Assume that all the fibers of $q$ are quadrics of corank $\leq 1$ and that the locus $Z\hookrightarrow B$ of critical values of $q$ is smooth.
\begin{itemize}
\item[(i)] When $d_q$ is even, we have the following equivalences of conjectures 
\begin{eqnarray*} \mathrm{C}(B) + \mathrm{C}(\widetilde{B}) \Leftrightarrow \mathrm{C}(Q) & \text{with} & \mathrm{C} \in \{\mathrm{C}^+, \mathrm{D}, \mathrm{V}, \mathrm{B}, \mathrm{W}, \mathrm{T}^l\, (l\neq 2), \mathrm{T}^p, \mathrm{ST}, \mathrm{P}, \mathrm{S}\}\,,
\end{eqnarray*}
where $\widetilde{B}$ stands for the discriminant twofold cover of $B$ (ramified over $Z$).
\item[(ii)] When $d_q$ is odd, $p\neq 2$, $d\leq 1$, and $k$ is algebraically closed or a finite field, we have the following equivalences of conjectures:
\begin{eqnarray*}
\mathrm{C}(B) + \mathrm{C}(Z) \Leftrightarrow \mathrm{C}(Q) & \text{with} & \mathrm{C} \in \{\mathrm{C}^+, \mathrm{D}, \mathrm{V}, \mathrm{B}, \mathrm{W}, \mathrm{T}^l, \mathrm{T}^p, \mathrm{ST}, \mathrm{P}, \mathrm{S}\}\,.
\end{eqnarray*}
\item[(iii)] When $d_q$ is odd and $p\neq 2$, we have the following implication of conjectures 
$$\{\mathrm{S}(U_j)\} + \{\mathrm{S}(\widetilde{Z}_j)\} \Rightarrow \mathrm{S}(Q)\,,$$ 
where $U_j$ is an affine open subscheme of $B$ and $\widetilde{Z}_j$ is a certain Galois twofold cover of $Z_j:=Z\cap U_j$ induced by the restriction of $q$ to $Z_j$.
\end{itemize}
\end{theorem}
Roughly speaking, Theorem \ref{thm:app1-1} relates the celebrated conjectures for the total space $Q$ with the celebrated conjectures for the base $B$. Items (i)-(ii) were proved in \cite[Thm.~1.2]{Crelle}\cite[Thm.~1.1(i)]{NCSchur} in the case of the conjectures $\mathrm{V}$ and $\mathrm{S}$. The proof of the other cases is similar. Item (iii) was proved in \cite[Thm.~1.1(ii)]{NCSchur}.
\begin{corollary}[Low-dimensional bases]\label{cor:app1-1}
Let $Q$ be as in Theorem \ref{thm:app1-1}.
\begin{itemize}
\item[(i)] When $d_q$ is even and $d\leq 1$, the following conjectures hold:
\begin{eqnarray*}
\mathrm{C}(Q) & \mathrm{with} & \mathrm{C} \in \{\mathrm{C}^+, \mathrm{D}, \mathrm{V}, \mathrm{B}, \mathrm{W}, \mathrm{T}^l \,(l\neq 2), \mathrm{T}^p, \mathrm{ST}, \mathrm{P}, \mathrm{S}\}\,.
\end{eqnarray*}
Moreover, $\mathrm{C}^+(Q)$ holds when $d\leq 2$, $\mathrm{D}(Q)$ holds when $d\leq 2$ or when $d\leq 4$ and $p=0$, and $\mathrm{V}(Q)$ holds when $d\leq 2$.
\item[(ii)]  When $d_q$ is odd, $p\neq 2$, $d\leq 1$, and $k$ is algebraically closed or a finite field, the following conjectures hold:
\begin{eqnarray*}
\mathrm{C}(Q)& \mathrm{with} & \mathrm{C} \in \{\mathrm{C}^+, \mathrm{D}, \mathrm{V}, \mathrm{B}, \mathrm{W}, \mathrm{T}^l, \mathrm{T}^p, \mathrm{ST}, \mathrm{P}, \mathrm{S}\}\,.
\end{eqnarray*}
\item[(iii)] When $d_q$ is odd, $p\neq 2$, and $d\leq 1$, the conjecture $\mathrm{S}(Q)$ holds. Moreover, when $d\leq 2$, we have the implication of conjectures $\mathrm{S}(B) \Rightarrow \mathrm{S}(Q)$.
\end{itemize}
\end{corollary}
\begin{remark}[Kimura-finiteness conjecture]
Assume that $B$ is a smooth $k$-curve.
\begin{itemize}
\item[(i)] When $d_q$ is even, it follows from Corollary \ref{cor:app1-1}(i) that the mixed motive $M(Q)_\bbQ$ is Schur-finite. As proved in \cite[Thm.~1.1(i)]{Kimura-curve}, $M(Q)_\bbQ$ is moreover Kimura-finite and $\mathrm{kim}(M(Q)_\bbQ)= d_q\cdot \mathrm{kim}(M(B)_\bbQ)+ \mathrm{kim}(M(\widetilde{B})_\bbQ)$.
\item[(ii)] When $d_q$ is odd, $p\neq 2$, and $k$ is algebraically closed or a finite field, it follows from Corollary \ref{cor:app1-1}(ii) that the mixed motive $M(Q)_\bbQ$ is Schur-finite. As proved in \cite[Thm.~1.1(ii)]{Kimura-curve}, the mixed motive $M(Q)_\bbQ$ is moreover Kimura-finite and $\mathrm{kim}(M(Q)_\bbQ)= (d_q+1)\cdot\mathrm{kim}(M(B)_\bbQ) + \#Z$.
\end{itemize}
\end{remark}
\begin{remark}[Bass-finiteness conjecture]
Let $k=\bbF_q$ be a finite field and $X$ a smooth $k$-scheme of finite type. In the seventies, Bass \cite{Bass} conjectured that the algebraic $K$-theory groups $K_n(X), n \geq 0$, are finitely generated. In the same vein, we can consider the mod $2$-torsion Bass-finiteness, where $K_n(X)$ is replaced by $K_n(X)_{1/2}$. As proved in \cite[Thm.~1.10]{NCSchur}, the above Theorem \ref{thm:app1-1} (items (i)-(iii)) holds similarly for the mod $2$-torsion Bass-finiteness conjecture.
\end{remark}
%-------------------------------------------------------------------------------
\subsection{Intersections of quadrics}
%-------------------------------------------------------------------------------
Let $k$ be a perfect base field of characteristic $p\geq 0$ and $X$ a smooth complete intersection of $m$ quadric hypersurfaces in $\bbP^n$. The linear span of these quadric hypersurfaces give rise to a flat quadric fibration $q\colon Q \to \bbP^{m-1}$ of relative dimension $n-1$.
\begin{theorem}\label{thm:app1-22}
Assume that all the fibers of $q$ are quadrics of corank $\leq 1$ and that the locus $Z \hookrightarrow \bbP^{m-1}$ of critical values of $q$ is smooth. Under these assumptions, we have the following equivalences/implications of conjectures:
\begin{eqnarray*}
\begin{cases} \mathrm{C}(Q) \Leftrightarrow \mathrm{C}(X) & 2m \leq n-1 \\
\mathrm{C}(Q) \Rightarrow \mathrm{C}(X) & 2m > n-1
\end{cases} & \mathrm{with} & \mathrm{C} \in \{\mathrm{C}^+, \mathrm{D}, \mathrm{V}, \mathrm{B}, \mathrm{W}, \mathrm{T}^l, \mathrm{T}^p, \mathrm{ST}, \mathrm{P}, \mathrm{S}\}\,.
\end{eqnarray*}
\end{theorem}
Intuitively speaking, Theorem \ref{thm:app1-22} shows that in order to solve the celebrated conjectures for intersections of quadrics, it suffices to solve the celebrated conjectures from quadric fibrations (and vice-versa).  This result was proved in \cite[Thm.~1.5]{NCSchur} in the case of the conjecture $\mathrm{S}$. The proof of the other cases is similar.
\begin{corollary}[Intersections of up-to-five quadrics]
Let $X$~be~as~in~Theorem~\ref{thm:app1-22}.
\begin{itemize}
\item[(i)] When $n$ is odd and $m\leq 2$, the following conjectures hold:
\begin{eqnarray*}
\mathrm{C}(X) & \mathrm{with} & \mathrm{C} \in \{\mathrm{C}^+, \mathrm{D}, \mathrm{V}, \mathrm{B}, \mathrm{W}, \mathrm{T}^l\,(l\neq 2), \mathrm{T}^p, \mathrm{ST}, \mathrm{P}, \mathrm{S}\}\,.
\end{eqnarray*}
Moreover, $\mathrm{C}^+(X)$ holds when $m\leq 3$, $\mathrm{D}(X)$ holds when $m\leq 3$ or when $m\leq 5$ and $p=0$, and $\mathrm{V}(X)$ holds when $m\leq 3$.
\item[(ii)] When $n$ is even, $p\neq 2$, $m\leq 2$, and $k$ is algebraically closed or a finite field, the following conjectures hold:
\begin{eqnarray*}
\mathrm{C}(X) & \mathrm{with} & \mathrm{C} \in \{\mathrm{C}^+, \mathrm{D}, \mathrm{V}, \mathrm{B}, \mathrm{W}, \mathrm{T}^l, \mathrm{T}^p, \mathrm{ST}, \mathrm{P}, \mathrm{S}\}\,.
\end{eqnarray*}
\item[(iii)] When $n$ is even, $p\neq 2$, and $m\leq 3$, the conjecture $\mathrm{S}(X)$ holds.
\end{itemize}
\end{corollary}
%-------------------------------------------------------------------------------
\subsection{Families of sextic du Val del Pezzo surfaces}
%-------------------------------------------------------------------------------
Let $k$ be a perfect base field of characteristic $p \geq 0$, $B$ a smooth proper $k$-scheme of dimension $d$ (in the case of conjecture $\mathrm{S}$ we assume solely that $B$ is smooth), and $f\colon X \to B$ a {\em family of sextic du Val del Pezzo surfaces}, i.e., a flat morphism such that for every geometric point $b \in B$ the associated fiber $X_b$ is a sextic du Val del Pezzo surface\footnote{Recall that a {\em sextic du Val del Pezzo surface} is a projective surface $S$ with at worst du Val singularities and whose ample anticanonical class $K_S$ is such that $K_S^2=6$.}. Following Kuznetsov \cite[\S5]{delPezzo}, let $\cM_2$, resp. $\cM_3$, be the relative moduli stack of semistable sheaves on fibers of $X$ over $B$ with Hilbert polynomial $h_2(t):=(3t+2)(t+1)$, resp. $h_3(t):=(3t+3)(t+1)$, and $Z_2$, resp. $Z_3$, the coarse moduli space~of~$\cM_2$,~resp.~$\cM_3$. %By construction, we have a finite flat morphism $Z_2 \to B$, resp. $Z_3 \to B$, of degree $3$, resp. $2$.
\begin{theorem}\label{thm:app1-2}
Assume that $p\neq 2,3$ and that $X$ is smooth. Under these assumptions, we have the following equivalences of conjectures:
\begin{eqnarray*}
\mathrm{C}(B) + \mathrm{C}(Z_2) + \mathrm{C}(Z_3) \Leftrightarrow \mathrm{C}(X) &\text{with} & \mathrm{C} \in \{\mathrm{C}^+, \mathrm{D}, \mathrm{V}, \mathrm{B}, \mathrm{W}, \mathrm{T}^l, \mathrm{T}^p, \mathrm{ST}, \mathrm{P}, \mathrm{S}\}\,.
\end{eqnarray*}
\end{theorem}
Roughly speaking, Theorem \ref{thm:app1-2} relates the celebrated conjectures for the total space $X$ with the celebrated conjectures for the base $B$. Theorem \ref{thm:app1-2} was proved in \cite[Thm.~1.7]{NCSchur} in the case of the conjecture $\mathrm{S}$. The proof of the other~cases~is~similar.
\begin{corollary}[Low-dimensional bases]
Let $X$ be as in Theorem \ref{thm:app1-2}. When $d\leq 1$, the following conjectures hold:
\begin{eqnarray*}
\mathrm{C}(X) & \mathrm{with} & \mathrm{C} \in \{\mathrm{C}^+, \mathrm{D}, \mathrm{V}, \mathrm{B}, \mathrm{W}, \mathrm{T}^l, \mathrm{T}^p, \mathrm{ST}, \mathrm{P}, \mathrm{S}\}\,.
\end{eqnarray*}
Moreover, $\mathrm{C}^+(X)$ holds when $d\leq 2$, $\mathrm{D}(X)$ holds when $d\leq 2$ or when $d\leq 4$ and $p=0$, and $\mathrm{V}(X)$ holds when $d\leq 2$. 
\end{corollary}

%-------------------------------------------------------------------------------
\subsection{Linear sections of Grassmannians}
%-------------------------------------------------------------------------------
Let $k$ be a base field of characteristic $p=0$, $W$ a $k$-vector space of dimension $6$ or $7$, $X:=\mathrm{Gr}(2,W)$ the Grassmannian variety of $2$-dimensional subspaces equipped with the Pl\"ucker embedding $\mathrm{Gr}(2,W) \hookrightarrow \bbP(\wedge^2(W))$, and $Y$ the Pfaffian variety $\mathrm{Pf}(4, W^\ast) \subset \bbP(\wedge^2(W^\ast))$. Given a linear subspace $L \subset \wedge^2 (W^\ast)$, consider the associated linear sections 
\begin{eqnarray*}
X_L:=X\times_{\bbP(\wedge^2(W))}\bbP(L^\perp)&& Y_L:=Y\times_{\bbP(\wedge^2(W^\ast))}\bbP(L)\,,
\end{eqnarray*}
where $L^\perp$ stands for the kernel of the induced homomorphism $\wedge^2(W) \twoheadrightarrow L^\ast$.
\begin{theorem}\label{thm:app1-3}
Assume that $X_L$ and $Y_L$ are smooth\footnote{The linear section $X_L$ is smooth if and only if the linear section $Y_L$ is smooth.}, and that $\mathrm{codim}(X_L)=\mathrm{dim}(L)$ and $\mathrm{codim}(Y_L)=\mathrm{dim}(L^\perp)$. Under these assumptions (which hold for a generic choice of $L$), we have the following equivalences of conjectures:
\begin{eqnarray*}
\mathrm{C}(X_L) \Leftrightarrow \mathrm{C}(Y_L) & \text{with} & \mathrm{C} \in \{\mathrm{C}^+, \mathrm{D}, \mathrm{V}, \mathrm{S}\}\,.
\end{eqnarray*}
\end{theorem}
Intuitively speaking, Theorem \ref{thm:app1-3} shows that in order to solve the celebrated conjectures for the linear section $X_L$, it suffices to solve the celebrated conjectures for the linear section $Y_L$ (and vice-versa). Theorem \ref{thm:app1-3} was proved in \cite[Thm.~1.7]{Crelle} in the case of the conjecture $\mathrm{V}$. The proof of the other cases is similar. 

When $\mathrm{dim}(W)=6$, we have $\mathrm{dim}(X_L)=8 - \mathrm{dim}(L)$ and $\mathrm{dim}(Y_L)=\mathrm{dim}(L)-2$. Moreover, in the case where $\mathrm{dim}(L)=5$, resp. $\mathrm{dim}(L)=6$, $X_L$ is a Fano threefold, resp. $K3$-surface, and $Y_L$ is a cubic threefold, resp. cubic fourfold.

When $\mathrm{dim}(W)=7$, we have $\mathrm{dim}(X_L)=10 - \mathrm{dim}(L)$ and $\mathrm{dim}(Y_L)=\mathrm{dim}(L)-4$. Moreover, in the case where $\mathrm{dim}(L)=5$, resp. $\mathrm{dim}(L)=6$, $X_L$ is a Fano fivefold, resp. Fano fourfold, and $Y_L$ is a curve of degree $42$, resp. surface of degree $42$. Furthermore, in the case where $\mathrm{dim}(L)=7$, $X_L$ and $Y_L$ are derived equivalent Calabi-Yau threefolds. In this particular latter case, Theorem \ref{thm:app1-3} follows then from the above Corollary \ref{cor:derived}.
\begin{corollary}[Low-dimensional sections]\label{cor:low}
Let $X_L$ be as in Theorem \ref{thm:app1-3}.
\begin{itemize}
\item[(i)] When $\mathrm{dim}(W)=6$ and $\mathrm{dim}(L)\leq 3$, the following conjectures $\mathrm{C}(X_L)$, with $\mathrm{C} \in \{\mathrm{C}^+, \mathrm{D}, \mathrm{V}, \mathrm{S}\}$, hold. Moreover, $\mathrm{C}^+(X_L)$ holds when $\mathrm{dim}(L)\leq 4$, $\mathrm{D}(X_L)$ holds when $\mathrm{dim}(L)\leq 6$, and $\mathrm{V}(X_L)$ holds when $\mathrm{dim}(L)\leq 4$.\
\item[(ii)] When $\mathrm{dim}(W)=7$ and $\mathrm{dim}(L)\leq 5$, the following conjectures $\mathrm{C}(X_L)$, with $\mathrm{C} \in \{\mathrm{C}^+, \mathrm{D}, \mathrm{V}, \mathrm{S}\}$, hold. Moreover, $\mathrm{C}^+(X_L)$ holds when $\mathrm{dim}(L)\leq 6$, $\mathrm{D}(X_L)$ holds when $\mathrm{dim}(L)\leq 8$, and $\mathrm{V}(X_L)$ holds when $\mathrm{dim}(L)\leq 6$.
\end{itemize}
\end{corollary}
\begin{remark}[Kimura-finiteness conjecture]
Let $X_L$ be as in Theorem \ref{thm:app1-3}. When $\mathrm{dim}(W)=6$ and $\mathrm{dim}(L)\leq 3$ or when $\mathrm{dim}(W)=7$ and $\mathrm{dim}(L)\leq 5$, it follows from Corollary \ref{cor:low} that the Chow motive $\mathfrak{h}(X_L)_\bbQ$ is Schur-finite. As proved in \cite[Thm.~1.10]{MRL}, the Chow motive $\mathfrak{h}(X_L)_\bbQ$ is moreover Kimura-finite.
\end{remark}
%-------------------------------------------------------------------------------
\subsection{Linear sections of Lagrangian Grassmannians}
%-------------------------------------------------------------------------------
Let $k$ be a base field of characteristic $p=0$, $W$ a $k$-vector space of dimension $6$ equipped with a symplectic form $\omega$, and $X:=\mathrm{LGr}(3,W)$ the associated Lagrangian Grassmannian of $3$-dimensional subspaces. The natural representation of the symplectic group $\mathrm{Sp}(\omega)$ on $\wedge^3(W)$ decomposes into a direct sum $W\oplus V$. Moreover, the classical Pl\"ucker embedding $\mathrm{Gr}(3,W) \hookrightarrow \bbP(\wedge^3(W))$ restricts to an embedding $\mathrm{LGr}(3,W) \hookrightarrow \bbP(V)$ of the Lagrangian Grassmannian. Consider also the classical projective dual variety $\mathrm{LGr}(3,6)^\vee\subseteq \bbP(V^\ast)$. This is a quartic hypersurface which is singular along a closed subvariety $Z$ of dimension $9$. Let us denote by $Y$ the open dense subset $\mathrm{LGr}(3,6)^\vee\backslash Z$. Given a linear subspace $L \subseteq V^\ast$ such that $\bbP(L)\cap Z = \emptyset$, consider the associated smooth linear sections $X_L:=X\times_{\bbP(V)}\bbP(L^\perp)$ and $Y_L:=Y \times_{\bbP(V^\ast)}\bbP(L)$.
\begin{theorem}\label{thm:app1-?}
Assume that $\mathrm{codim}(X_L)=\mathrm{dim}(L)$ and $\mathrm{codim}(Y_L)=\mathrm{dim}(L^\perp)$. Under these assumption (which hold for a generic choice of $L$), we have the following equivalences of conjectures:
\begin{eqnarray*}
\mathrm{C}(X_L) \Leftrightarrow \mathrm{C}(Y_L) &\text{with} & \mathrm{C} \in \{\mathrm{C}^+, \mathrm{D}, \mathrm{V}, \mathrm{S}\}\,.
\end{eqnarray*}  
\end{theorem}
Intuitively speaking, Theorem \ref{thm:app1-?} shows that in order to solve the celebrated conjectures for the linear section $X_L$, it suffices to solve the celebrated conjectures for the linear section $Y_L$ (and vice-versa).
Theorem \ref{thm:app1-?} was proved in \cite[Thm.~1.5]{MRL} in the case of the conjecture $\mathrm{S}$. The proof of the other cases is similar. 

We have $\mathrm{dim}(X_L)=6 - \mathrm{dim}(L)$ and $\mathrm{dim}(Y_L)=\mathrm{dim}(L)-2$. Moreover, in the case where $\mathrm{dim}(L)=3$, resp. $\mathrm{dim}(L)=4$, $X_L$ is a Fano threefold, resp. $K3$-surface of degree $16$, and $Y_L$ is a plane quartic, resp. $K3$-surface~of~degree~$4$.
\begin{corollary}[Low-dimensional sections]\label{cor:low1}
Let $X_L$ be as in Theorem \ref{thm:app1-?}. When $\mathrm{dim}(L)\leq 3$, the following conjectures $\mathrm{C}(X_L)$, with $ \mathrm{C} \in \{\mathrm{C}^+, \mathrm{D}, \mathrm{V}, \mathrm{S}\}$, hold. Moreover, $\mathrm{C}^+(X_L)$ holds when $\mathrm{dim}(L)\leq 4$, $\mathrm{D}(X_L)$ holds when $\mathrm{dim}(L)\leq 6$, and $\mathrm{V}(X_L)$ holds when $\mathrm{dim}(L)\leq 4$.
\end{corollary}
\begin{remark}[Kimura-finiteness conjecture]
Let $X_L$ be as in Theorem \ref{thm:app1-?}. When $\mathrm{dim}(L)\leq 3$, it follows from Corollary \ref{cor:low1} that the Chow motive $\mathfrak{h}(X_L)_\bbQ$ is Schur-finite. As proved in \cite[Thm.~1.10]{MRL}, $\mathfrak{h}(X_L)_\bbQ$ is moreover Kimura-finite.
\end{remark}
%-------------------------------------------------------------------------------
\subsection{Linear sections of spinor varieties}
%-------------------------------------------------------------------------------
Let $k$ be a base field of characteristic $p=0$, $W$ a $k$-vector space of dimension $10$ equipped with a nondegenerate quadratic form $q \in \mathrm{Sym}^2(W^\ast)$, and $X:=\mathrm{OGr}_+(5,W)$ and $Y:=\mathrm{OGr}_-(5,W)$ the connected components of the orthogonal Grassmannian of $5$-dimensional subspaces. These are called the {\em spinor varieties}. By construction, we have a canonical embedding $\mathrm{OGr}_+(5,W) \hookrightarrow \bbP(V)$, where $V$ stands for the corresponding half-spinor representation of the spin-group $\mathrm{Spin}(W)$. In the same vein, making use of the isomorphism $\bbP(V) \simeq \bbP(V^\ast)$ induced by the nondegenerate quadratic form $q$, we have the embedding $\mathrm{OGr}_-(5,W) \hookrightarrow \bbP(V^\ast)$. Given a linear subspace $L \subseteq V^\ast$, consider the associated linear sections $X_L:=X\times_{\bbP(V)}\bbP(L^\perp)$ and $Y_L:=Y \times_{\bbP(V^\ast)}\bbP(L)$.
\begin{theorem}\label{thm:app1-??}
Assume that $X_L$ and $Y_L$ are smooth, and that $\mathrm{codim}(X_L)=\mathrm{dim}(L)$ and $\mathrm{codim}(Y_L)=\mathrm{dim}(L^\perp)$. Under these assumptions (which hold for a generic choice of $L$), we have the following equivalences of conjectures:
 \begin{eqnarray*}
\mathrm{C}(X_L) \Leftrightarrow \mathrm{C}(Y_L) &\text{with} & \mathrm{C} \in \{\mathrm{C}^+, \mathrm{D}, \mathrm{V}, \mathrm{S}\}\,.
\end{eqnarray*} 
\end{theorem}
Intuitively speaking, Theorem \ref{thm:app1-?} shows that in order to solve the celebrated conjectures for the linear section $X_L$, it suffices to solve the celebrated conjectures for the linear section $Y_L$ (and vice-versa).
Theorem \ref{thm:app1-??} was proved in \cite[Thm.~1.6]{MRL} in the case of the conjecture $\mathrm{S}$. The proof of the other cases is similar. 

We have $\mathrm{dim}(X_L)=10-\mathrm{dim}(L)$ and $\mathrm{dim}(Y_L)=\mathrm{dim}(L)-6$. Moreover, in the case where $\mathrm{dim}(L)=7$, $X_L$ is a Fano threefold and $Y_L$ is a curve of genus $7$. Furthermore, in the case where $\mathrm{dim}(L)=8$, $X_L$ and $Y_L$ are derived equivalent $K3$-surfaces of degree $12$. In this particular latter case, Theorem \ref{thm:app1-??} follows then from the above Corollary \ref{cor:derived}.
\begin{corollary}[Low-dimensional sections]\label{cor:low2}
Let $X_L$ be as in Theorem \ref{thm:app1-??}. When $\mathrm{dim}(L)\leq 7$, the following conjectures $\mathrm{C}(X_L)$, with $\mathrm{C} \in \{\mathrm{C}^+, \mathrm{D}, \mathrm{V}, \mathrm{S}\}$, hold. Moreover, $\mathrm{C}^+(X_L)$ holds when $\mathrm{dim}(L)\leq 8$, $\mathrm{D}(X_L)$ holds when $\mathrm{dim}(L)\leq 10$, and $\mathrm{V}(X_L)$ holds when $\mathrm{dim}(L)\leq 8$.
\end{corollary}
\begin{remark}[Kimura-finiteness conjecture]
Let $X_L$ be as in Theorem \ref{thm:app1-??}. When $\mathrm{dim}(L)\leq 7$, it follows from Corollary \ref{cor:low2} that the Chow motive $\mathfrak{h}(X_L)_\bbQ$ is Schur-finite. As proved in \cite[Thm.~1.10]{MRL}, $\mathfrak{h}(X_L)_\bbQ$ is moreover Kimura-finite.
\end{remark}
%-------------------------------------------------------------------------------
\subsection{Linear sections of determinantal varieties}
%-------------------------------------------------------------------------------
Let $k$ be a perfect base field of characteristic $p\geq 0$, $U_1$ and $U_2$ two finite-dimensional $k$-vector spaces of dimensions $d_1$ and $d_2$, respectively, $V:=U_1 \otimes U_2$, and $0 < r < d_1$ an integer. Consider the determinantal variety $\cZ^r_{d_1,d_2}\subset \bbP(V)$ defined as the locus of those matrices $U_2 \to U_1^\ast$ with rank $\leq r$; recall that this condition can be described as the vanishing of the $(r+1)$-minors of the matrix of indeterminates:
$$ \begin{pmatrix}
x_{1,1} & \cdots & x_{1,d_2} \\
\vdots & \ddots & \vdots \\
x_{d_1, 1} & \cdots & x_{d_1, d_2}
\end{pmatrix}\,.
$$
\begin{example}[Segre varieties]%\label{ex:Segre1}
In the particular case where $r=1$, the determinantal varieties reduce to the classical Segre varieties. Concretely, $\cZ_{d_1, d_2}^1$ reduces to the image of Segre homomorphism $\bbP(U_1) \times \bbP(U_2) \to \bbP(V)$. For example, $\cZ_{2,2}^1$ is the classical quadric hypersurface:
$$\{[x_{1,1}:x_{1,2}: x_{2,1}: x_{2,2}]\,|\,\mathrm{det}\begin{pmatrix} x_{1,1} & x_{1,2} \\ x_{2,1} & x_{2,2} \end{pmatrix}=0 \}\subset \bbP^3\,.$$
\end{example}
In contrast with the Segre varieties, the determinantal varieties $\cZ^r_{d_1, d_2}$, with $r\geq 2$, are not smooth. The singular locus of $\cZ^r_{d_1, d_2}$ consists of those matrices $U_2 \to U_1^\ast$ with rank $<r$, \ie it agrees with the closed subvariety $\cZ^{r-1}_{d_1, d_2}$. Nevertheless, it is well-known that $\cZ^r_{d_1, d_2}$ admits a canonical Springer resolution of singularities $X:=\cX_{d_1, d_2}^r \to \cZ^r_{d_1, d_2}$. Dually, consider the variety $\cW^r_{d_1, d_2}\subset \bbP(V^\ast)$, defined as the locus of those matrices $U^\ast_2 \to U_1$ with corank $\geq r$, and the associated canonical Springer resolution of singularities $Y:=\cY^r_{d_1, d_2} \to \cW^r_{d_1, d_2}$. Given a linear subspace $L\subseteq V^\ast$, consider the associated linear sections $X_L:=X\times_{\bbP(V)}\bbP(L^\perp)$ and $Y_L:=Y \times_{\bbP(V^\ast)}\bbP(L)$. Note that whenever $\bbP(L^\perp)$ does not intersects the singular locus of $\cZ^r_{d_1, d_2}$, we have $X_L=\bbP(L^\perp) \cap \cZ^r_{d_1, d_2}$, i.e., $X_L$ is a linear section of the determinantal variety $\cZ^r_{d_1, d_2}$.
\begin{theorem}\label{thm:app1-4}
Assume that $X_L$ and $Y_L$ are smooth, and that $\mathrm{codim}(X_L)=\mathrm{dim}(L)$ and $\mathrm{codim}(Y_L)=\mathrm{dim}(L^\perp)$. Under these assumptions (which hold for a generic choice of $L$), we have the following equivalences of conjectures:
\begin{eqnarray*}
\mathrm{C}(X_L) \Leftrightarrow \mathrm{C}(Y_L) & \text{with} & \mathrm{C} \in \{\mathrm{C}^+, \mathrm{D}, \mathrm{V}, \mathrm{B}, \mathrm{W}, \mathrm{T}^l, \mathrm{T}^p, \mathrm{ST}, \mathrm{P}, \mathrm{S}\}\,.
\end{eqnarray*}
\end{theorem}
Intuitively speaking, Theorem \ref{thm:app1-?} shows that in order to solve the celebrated conjectures for the linear section $X_L$, it suffices to solve the celebrated conjectures for the linear section $Y_L$ (and vice-versa).
Theorem \ref{thm:app1-4} was proved in \cite[Cor.~2.4]{NCpositive-CD}, resp. in \cite[Cor.~1.7]{NCTate}, in the case of the conjectures $\mathrm{C}^+$ and $\mathrm{D}$, resp. in the case of the conjectures $\mathrm{B}$, $\mathrm{T}^l$, $\mathrm{T}^p$, and $\mathrm{P}$. The proof of the other cases is similar. 

By construction, we have the following equalities: 
$$\mathrm{dim}(X_L)= r(d_1 + d_2 - r)-1-\mathrm{dim}(L)\quad \mathrm{dim}(Y_L)= r(d_1-d_2 -r)-1 + \mathrm{dim}(L)\,.$$
Moreover, in the case where $\mathrm{dim}(L)=d_2r$, the linear sections $X_L$ and $Y_L$ are derived equivalent Calabi-Yau varieties. In this particular latter case, Theorem \ref{thm:app1-4} follows then from the above Corollary \ref{cor:derived}.
\begin{corollary}[High-dimensional sections]\label{cor:determinantal}
Let $X_L$ be as in Theorem \ref{thm:app1-4}. When $\mathrm{dim}(L)\leq 2 - r(d_1 - d_2 - r)$, the following conjectures hold:
\begin{eqnarray}
\mathrm{C}(X_L) & \mathrm{with} & \mathrm{C} \in \{\mathrm{C}^+, \mathrm{D}, \mathrm{V}, \mathrm{B}, \mathrm{W}, \mathrm{T}^l, \mathrm{T}^p, \mathrm{ST}, \mathrm{P}, \mathrm{S}\}\,.
\end{eqnarray}
Moreover, $\mathrm{C}^+(X_L)$ holds when $\mathrm{dim}(L)\leq 3 - r(d_1 - d_2 - r)$, $\mathrm{D}(X_L)$ holds when $\mathrm{dim}(L)\leq 3 - r(d_1 - d_2 - r)$ or when $\mathrm{dim}(L)\leq 5 - r(d_1 - d_2 - r)$ and $p=0$, and $\mathrm{V}(X_L)$ holds when $\mathrm{dim}(L)\leq 3 - r(d_1 - d_2 - r)$.
\end{corollary}
\begin{example}[Segre varieties]\label{ex:Segre}
Let $r=1$. Thanks to Corollary \ref{cor:determinantal}, when $\mathrm{dim}(L)=3 - d_1 + d_2$, the following conjectures hold:
\begin{eqnarray*}
\mathrm{C}(X_L) & \mathrm{with} & \mathrm{C}\in \{\mathrm{C}^+, \mathrm{D}, \mathrm{V}, \mathrm{B}, \mathrm{W}, \mathrm{T}^l, \mathrm{T}^p, \mathrm{ST}, \mathrm{P}, \mathrm{S}\}\,.
\end{eqnarray*}
In all these cases, $X_L$ is a smooth linear section of the Segree variety $\cZ^1_{d_1, d_2}$. Moreover, $X_L$ is Fano if and only if $\mathrm{dim}(L)<d_1$. Furthermore, $\mathrm{dim}(X_L)=2d_1-5$. Therefore, by letting $d_1 \to \infty$ and by keeping $\mathrm{dim}(L)$ fixed, we obtain infinitely many examples of smooth proper $k$-schemes $X_L$, of arbitrary high dimension, satisfying the celebrated conjectures of \S\ref{sec:conjectures}. Note that in the particular case of the Weil conjecture (and in the case of the Grothendieck standard conjecture of type $\mathrm{C}^+$ over a finite field) this proof avoids all the technical tools used by Deligne.
\end{example}
\begin{example}[Rational normal scrolls]
Let $r=1$, $d_1=4$ and $d_2=2$. In this particular case, the Segre variety $\cZ^1_{4,2}\subset \bbP^7$ agrees with the rational normal $4$-fold scroll $S_{1, 1, 1, 1}$. Choose a linear subspace $L\subseteq V^\ast$ of dimension $1$ such that the hyperplane $\bbP(L^\perp) \subset \bbP^7$ does not contains any $3$-plane of the rulling of $S_{1, 1, 1, 1}$; this condition holds for a generic choice of $L$. In this case, the linear section $X_L$ agrees with the $3$-fold scroll $S_{1, 1, 2}$. Hence, thanks to Example \ref{ex:Segre}, we conclude that the following conjectures hold:
\begin{eqnarray*}
\mathrm{C}(S_{1, 1, 2}) & \mathrm{with} & \mathrm{C}\in \{\mathrm{C}^+, \mathrm{D}, \mathrm{V}, \mathrm{B}, \mathrm{W}, \mathrm{T}^l, \mathrm{T}^p, \mathrm{ST}, \mathrm{P}, \mathrm{S}\}\,.
\end{eqnarray*}
\end{example}
\begin{example}[Square matrices]\label{ex:square}
Let $d_1=d_2$. Thanks to Corollary \ref{cor:determinantal}, when $\mathrm{dim}(L)=2 + r^2$, the following conjectures hold:
\begin{eqnarray*}
\mathrm{C}(X_L) & \mathrm{with} & \mathrm{C}\in \{\mathrm{C}^+, \mathrm{D}, \mathrm{V}, \mathrm{B}, \mathrm{W}, \mathrm{T}^l, \mathrm{T}^p, \mathrm{ST}, \mathrm{P}, \mathrm{S}\}\,.
\end{eqnarray*}
In all these cases, we have $\mathrm{dim}(X_L)=2r(d_1-r)-3$. Therefore, by letting $d_1\to \infty$ and by keeping $\mathrm{dim}(L)$ fixed, we obtain infinitely many new examples of smooth proper $k$-schemes $X_L$, of arbitrary high dimension, satisfying the celebrated conjectures of \S\ref{sec:conjectures}. Similarly to Example \ref{ex:Segre}, note that in the particular case of the Weil conjecture (and in the case of the Grothendieck standard conjecture of type $\mathrm{C}^+$ over a finite field) this proof avoids all the technical tools used by Deligne.
\end{example}
\begin{remark}[Kimura-finiteness conjecture]
Let $X_L$ be as in Theorem \ref{thm:app1-4}. When $\mathrm{dim}(L)= 2 - r(d_1 - d_2 -r)$, it follows from Corollary \ref{cor:determinantal} that the Chow motive $\mathfrak{h}(X_L)_\bbQ$ is Schur-finite. A similar proof shows that $\mathfrak{h}(X_L)_\bbQ$ is Kimura-finite.
\end{remark}
%-------------------------------------------------------------------------------

%-------------------------------------------------------------------------------
\section{Applications to noncommutative geometry}\label{sec:applications2}
%-------------------------------------------------------------------------------
The theorems of \S\ref{sec:NCcounterparts} enabled also a proof of the noncommutative counterparts of the celebrated conjectures of Grothendieck, Voevodsky, Beilinson, Weil, Tate, Kimura, and Schur, in many interesting cases. In this section, we describe some of these interesting cases. Similarly to \S\ref{sec:applications1}, we will often use the letter $\mathrm{C}$ to denote one of the celebrated conjectures $\{\mathrm{C}^+, \mathrm{D}, \mathrm{V}, \mathrm{B},\mathrm{W}, \mathrm{T}^l, \mathrm{T}^p, \mathrm{ST}, \mathrm{P}, \mathrm{K}, \mathrm{S}\}$. Moreover, given a smooth (proper) algebraic stack $\cX$, we will write $\mathrm{C}_{\mathrm{nc}}(\cX)$ instead~of~$\mathrm{C}_{\mathrm{nc}}(\perf_\dg(\cX))$.
%-------------------------------------------------------------------------------
\subsection{Finite-dimensional algebras of finite global dimension}
%-------------------------------------------------------------------------------
Let $k$ be a perfect base field of characteristic $p\geq 0$ and $A$ a finite-dimensional $k$-algebra of finite global dimension. Examples include path algebras of finite quivers without oriented cycles as well as their quotients by admissible ideals.
\begin{theorem}\label{thm:app2-1}
The following conjectures hold:
\begin{eqnarray*}
\mathrm{C}_{\mathrm{nc}}(A) & \text{with} & \mathrm{C} \in \{\mathrm{C}^+, \mathrm{D}, \mathrm{V}, \mathrm{B}, \mathrm{W}, \mathrm{T}^p, \mathrm{ST}, \mathrm{P}, \mathrm{K}, \mathrm{S}\}\,.
\end{eqnarray*}
Moreover, when $k=\overline{k}$, the conjecture $\mathrm{T}^l_{\mathrm{nc}}(A)$ also holds.
\end{theorem}
Theorem \ref{thm:app2-1} was proved in \cite[Thm.~3.1]{NCWeil} in the case of the conjectures $\mathrm{W}$ and $\mathrm{ST}$. The proof of the other cases is similar.
%-------------------------------------------------------------------------------
\subsection{Semi-orthogonal decompositions}
%-------------------------------------------------------------------------------
Let $k$ be a base field of characteristic $p\geq 0$ and $\cB, \cC \subseteq \cA$ smooth proper $k$-linear dg categories inducing a semi-orthogonal decomposition $\dgHo(\cA)=\langle \dgHo(\cB), \dgHo(\cC)\rangle$ in the sense of Bondal-Orlov \cite{BO}.
\begin{theorem}\label{thm:semi}
We have the following equivalences of conjectures
\begin{eqnarray*}
\mathrm{C}_{\mathrm{nc}}(\cB) + \mathrm{C}_{\mathrm{nc}}(\cC) \Leftrightarrow \mathrm{C}_{\mathrm{nc}}(\cA) & \mathrm{with} & \mathrm{C} \in \{\mathrm{C}^+, \mathrm{D}, \mathrm{V}, \mathrm{B}, \mathrm{W}, \mathrm{T}^l, \mathrm{T}^p, \mathrm{ST}, \mathrm{P}, \mathrm{K}, \mathrm{S}\}\,.
\end{eqnarray*}
\end{theorem}
Intuitively speaking, Theorem \ref{thm:semi} shows that the noncommutative counterparts of the celebrated conjectures are additive with respect to semi-orthogonal decompositions. Theorem \ref{thm:semi} was proved in \cite[Thm.~3.2]{NCWeil} in the case of the conjectures $\mathrm{W}$ and $\mathrm{ST}$. The proof of the other cases is similar.%------------------------------------------------------------------------------
\subsection{Calabi-Yau dg categories associated to hypersurfaces}
%-------------------------------------------------------------------------------
Let $k$ be a base field of characteristic $p\geq 0$ and $X \subset \bbP^n$ a smooth hypersurface of degree $\mathrm{deg}(X) \leq n+1$. Following Kuznetsov \cite{Kuznetsov-CY}, we have a semi-orthogonal decomposition:
$$\perf(X) = \langle \cT(X), \cO_X, \ldots, \cO_X(n-\mathrm{deg}(X))\rangle\,.$$ 
Moreover, the associated $k$-linear dg category $\cT_\dg(X)$, defined as the dg enhancement of $\cT(X)$ induced from $\perf_\dg(X)$, is a smooth proper Calabi-Yau dg category\footnote{In the particular case where $n=5$ and $\mathrm{deg}(X)=3$, the dg categories $\cT_\dg(X)$ obtained in this way are usually called ``noncommutative $K3$-surfaces'' because they share many of the key properties of the dg categories of perfect complexes of the classical $K3$-surfaces.} of fractional dimension $\frac{(n+1)(\mathrm{deg}(X)-2)}{\mathrm{deg}(X)}$. By combining Theorem \ref{thm:semi} with the Theorems of \S\ref{sec:NCcounterparts}, we hence obtain the following result:
\begin{corollary}\label{thm:app2-2}
We have the equivalences of conjectures
\begin{eqnarray*}
\mathrm{C}(X)\Leftrightarrow \mathrm{C}_{\mathrm{nc}}(\cT_\dg(X)) & \text{with} & \mathrm{C} \in \{\mathrm{C}^+, \mathrm{D}, \mathrm{V}, \mathrm{B}, \mathrm{W}, \mathrm{T}^l, \mathrm{T}^p, \mathrm{ST}, \mathrm{P}, \mathrm{S}\}
\end{eqnarray*}
as well as the implication $\mathrm{K}(X) \Rightarrow \mathrm{K}_{\mathrm{nc}}(\cT_\dg(X))$.
\end{corollary}
Roughly speaking, Corollary \ref{thm:app2-2} shows that in what concerns the celebrated conjectures, there is no difference between the hypersurface $X$ and the associated Calabi-Yau dg category $\cT_{\mathrm{dg}}(X)$.
%-------------------------------------------------------------------------------
\subsection{Root stacks}
%-------------------------------------------------------------------------------
Let $k$ be a base field of characteristic $p\geq 0$, $X$ a smooth proper $k$-scheme of dimension $d$ (in the case of conjecture $\mathrm{S}$ we assume solely that $X$ is smooth), $\cL$ a line bundle on $X$, $\varsigma \in \Gamma(X,\cL)$ a global section, and $n\geq 1$ an integer. Following Cadman \cite[\S2.2]{Cadman}, the associated {\em root stack} is defined as the fiber-product
$$
\xymatrix{
\cX:=\sqrt[n]{(\cL,\varsigma)/X} \ar[d]_-f \ar[r]& [\bbA^1/\bbG_m] \ar[d]^-{\theta_n} \\
X \ar[r]_-{(\cL,\varsigma)} & [\bbA^1/\bbG_m]\,,
}
$$
where $\theta_n$ stands for the morphism induced by the $n^{\mathrm{th}}$ power map on $\bbA^1$ and $\bbG_m$. As proved by Ishii-Ueda in \cite[Thm.~1.6]{IU}, whenever the zero locus $Z \hookrightarrow X$ of $\varsigma$ is smooth, we have a semi-orthogonal decomposition 
$$\perf(\cX) = \langle \perf(Z)_{n-1}, \ldots, \perf(Z)_1, f^\ast(\perf(X))\rangle\,,$$
where all the categories $\perf(Z)_j$ are (Fourier-Mukai) equivalent to $\perf(Z)$. Hence, by combining Theorem \ref{thm:semi} with the Theorems of \S\ref{sec:NCcounterparts}, we obtain the following results:
\begin{corollary}\label{thm:app2-22}
Assume that the zero locus $Z \hookrightarrow X$ of the global section $\varsigma$ is smooth. Under this assumption, we have the equivalences of conjectures
\begin{eqnarray*}
\mathrm{C}(X) + \mathrm{C}(Z) \Leftrightarrow \mathrm{C}_{\mathrm{nc}}(\cX) & \text{with} & \mathrm{C} \in \{\mathrm{C}^+, \mathrm{D}, \mathrm{V}, \mathrm{B}, \mathrm{W},  \mathrm{T}^l, \mathrm{T}^p, \mathrm{ST}, \mathrm{P}, \mathrm{S}\} 
\end{eqnarray*}
as well as the implication $\mathrm{K}(X) + \mathrm{K}(Z) \Rightarrow \mathrm{K}_{\mathrm{nc}}(\cX)$.
\end{corollary}
\begin{corollary}[Low-dimensional root stacks]
Let $\cX$ be as in Corollary \ref{thm:app2-22}. When $d\leq 1$, the following conjectures hold:
\begin{eqnarray*}
\mathrm{C}_{\mathrm{nc}}(\cX) & \mathrm{with} & \mathrm{C} \in \{\mathrm{C}^+, \mathrm{D}, \mathrm{V}, \mathrm{B}, \mathrm{W},  \mathrm{T}^l, \mathrm{T}^p, \mathrm{ST}, \mathrm{P}, \mathrm{K}, \mathrm{S}\}\,.
\end{eqnarray*}
Moreover, $\mathrm{C}^+_{\mathrm{nc}}(\cX)$ holds when $d\leq 2$, $\mathrm{D}_{\mathrm{nc}}(\cX)$ holds when $d\leq 2$ or when $d\leq 4$ and $p=0$, and $V_{\mathrm{nc}}(\cX)$ holds when $d\leq 2$.
\end{corollary}
%-------------------------------------------------------------------------------
\subsection{Global orbifolds}\label{sub:global}
%-------------------------------------------------------------------------------
Let $k$ be a base field of characteristic $p\geq 0$, $G$ a finite group of order $n$, $X$ a smooth proper $k$-scheme of dimension $d$ equipped with a $G$-action (in the case of conjecture $\mathrm{S}$ we assume solely that $X$ is smooth), and $\cX:=[X/G]$ the associated global orbifold.
\begin{theorem}\label{thm:app2-3}
Assume that $p\nmid n$ and that $k$ contains the $n^{\mathrm{th}}$ roots of unity. Under these assumptions, we have the following implications of conjectures
\begin{eqnarray*}
\sum_{\sigma \subseteq G} \mathrm{C}(X^\sigma) \Rightarrow \mathrm{C}_{\mathrm{nc}}(\cX) & \mathrm{with} & \mathrm{C} \in \{\mathrm{C}^+, \mathrm{D}, \mathrm{V}, \mathrm{B}, \mathrm{W},  \mathrm{T}^l\,(l\nmid n), \mathrm{T}^p, \mathrm{ST}, \mathrm{P}, \mathrm{K}, \mathrm{S}\}\,,
\end{eqnarray*} 
where $\sigma$ is a cyclic subgroup of $G$.
\end{theorem}
Intuitively speaking, Theorem \ref{thm:app2-3} shows that in order to solve the noncommutative counterparts of the celebrated conjectures for the global orbifold $\cX$, it suffices to solve the celebrated conjectures for the underlying scheme $X$. Theorem \ref{thm:app2-3} was proved in \cite[Thm.~9.2]{Orbifold}\cite[Thm.~3.1]{NCpositive-CD}, resp. in \cite[Thm.~1.16]{NCTate}, in the case of the conjectures $\mathrm{C}^+$, $\mathrm{D}$, and $\mathrm{V}$, resp. in the case of the conjectures $\mathrm{B}$, $\mathrm{T}^l\, (l\nmid n)$, $\mathrm{T}^p$, and $\mathrm{P}$. The proof of the other cases is similar.
\begin{corollary}[Low-dimensional global orbifolds]
Let $\cX$ be as in Theorem \ref{thm:app2-3}. When $d\leq 1$, the following conjectures hold:
\begin{eqnarray*}
\mathrm{C}_{\mathrm{nc}}(\cX) & \mathrm{with} & \mathrm{C} \in \{\mathrm{C}^+, \mathrm{D}, \mathrm{V}, \mathrm{B}, \mathrm{W},  \mathrm{T}^l\,(l\nmid n), \mathrm{T}^p, \mathrm{ST}, \mathrm{P}, \mathrm{K}, \mathrm{S}\}\,.
\end{eqnarray*}
Moreover, $\mathrm{C}^+_{\mathrm{nc}}$ holds when $d\leq 2$, $\mathrm{D}_{\mathrm{nc}}(\cX)$ holds when $d\leq 2$ or when $d\leq 4$ and $p=0$, $\mathrm{V}_{\mathrm{nc}}(\cX)$ holds when $d\leq 2$ or when $X$ is an abelian $3$-fold and $p=0$, and $\mathrm{T}_{\mathrm{nc}}^l(\cX)$ (with $l\nmid n$) and $\mathrm{T}^p_{\mathrm{nc}}(\cX)$ hold when $X$ is a $K3$-surface.
\end{corollary}
\begin{corollary}[Abelian $G$-varieties]
Let $\cX$ be as in Theorem \ref{thm:app2-3}. When $X$ is an abelian variety and $G$ acts by group homomorphisms\footnote{For example, in the case where $G=\bbZ/2$, we can consider the canonical involution $a \mapsto -a$.}, the conjectures $\mathrm{C}_{\mathrm{nc}}(\cX)$, with $\mathrm{C} \in \{\mathrm{C}^+, \mathrm{D}, \mathrm{W}, \mathrm{K}, \mathrm{S}\}$, hold. Moreover, when $d \leq 3$, the conjectures $\mathrm{C}_{\mathrm{nc}}(\cX)$, with $ \mathrm{C} \in \{\mathrm{B}, \mathrm{T}^l\,(l\nmid n), \mathrm{T}^p, \mathrm{ST}, \mathrm{P}\}$, also hold.
\end{corollary}
%-------------------------------------------------------------------------------
\subsection{Twisted global orbifolds}
%-------------------------------------------------------------------------------
Let $k$, $G$, $X$ (of dimension $d$), and $\cX:=[X/G]$, be as in \S\ref{sub:global}. In this subsection we consider the case where the global orbifold $\cX$ is equipped with a sheaf of Azumaya algebras $\cF$ of rank $r$. In other words, $\cF$ is a $G$-equivariant sheaf of Azumaya algebras of rank $r$ over $X$. Similarly to the dg category $\perf_\dg(\cX)$, we can also consider the dg category $\perf_\dg(\cX;\cF)$ of perfect complexes of $\cF$-modules. In what follows, we will write $\mathrm{C}_{\mathrm{nc}}(\cX;\cF)$ instead of $\mathrm{C}_{\mathrm{nc}}(\perf_\dg(\cX;\cF))$. The next result is the ``twisted'' version of the Theorem \ref{thm:app2-3}:
\begin{theorem}\label{thm:app2-4}
Assume that $p \nmid nr$ and that $k$ contains the $n^{\mathrm{th}}$ roots of unity. Under these assumptions, we have the following implications of conjectures
\begin{eqnarray*}
\sum_{\sigma \subseteq G} \mathrm{C}(Y_\sigma) \Rightarrow \mathrm{C}_{\mathrm{nc}}(\cX;\cF) & \text{with} & \mathrm{C} \in \{\mathrm{C}^+, \mathrm{D}, \mathrm{V}, \mathrm{B}, \mathrm{W}, \mathrm{T}^l\,(l\nmid nr), \mathrm{T}^p, \mathrm{ST}, \mathrm{P}, \mathrm{K}, \mathrm{S}\}\,,
\end{eqnarray*} 
where $\sigma$ is a cyclic subgroup of $G$ and $Y_\sigma$ is a certain $\sigma^\vee$-Galois cover of $X^\sigma$ induced by the restriction of $\cF$ to $X^\sigma$.
\end{theorem}
Intuitively speaking, Theorem \ref{thm:app2-4} shows that in order to solve the noncommutative counterparts of the celebrated conjectures for the twisted global orbifold $(\cX;\cF)$, it suffices to solve the celebrated conjectures for certain Galois covers of the underlying scheme $X$.
Theorem \ref{thm:app2-4} was proved in \cite[Thm.~1.23]{NCTate} in the case of the conjectures $\mathrm{B}$, $\mathrm{T}^l\, (l\nmid nr)$, $\mathrm{T}^p$, and $\mathrm{P}$. The proof of the other cases is similar.

\begin{corollary}[Low-dimensional twisted global orbifolds]
Let $\cX$ and $\cF$ be as in Theorem \ref{thm:app2-4}. When $d\leq 1$, the following conjectures hold:
\begin{eqnarray*}
\mathrm{C}_{\mathrm{nc}}(\cX;\cF) & \mathrm{with} & \mathrm{C} \in \{\mathrm{C}^+, \mathrm{D}, \mathrm{V}, \mathrm{B}, \mathrm{W}, \mathrm{T}^l\,(l\nmid nr), \mathrm{T}^p, \mathrm{ST}, \mathrm{P}, \mathrm{K}, \mathrm{S}\}\,.
\end{eqnarray*}
Moreover, $\mathrm{C}_{\mathrm{nc}}^+(\cX;\cF)$ holds when $d \leq 2$, $\mathrm{D}_{\mathrm{nc}}(\cX;\cF)$ holds when $d\leq 2$ or when $d\leq 4$ and $p=0$, and $\mathrm{V}_{\mathrm{nc}}(\cX;\cF)$ holds when $d\leq 2$.
\end{corollary}

%-------------------------------------------------------------------------------
\subsection{Intersections of bilinear divisors}
%-------------------------------------------------------------------------------
Let $k$ be a base field of characteristic $p\geq 0$ and $V$ a finite-dimensional $k$-vector space. Consider the canonical action of $\bbZ/2$ on $\bbP(V) \times \bbP(V)$ and the associated global orbifold $\cX:=[(\bbP(V) \times \bbP(V))/(\bbZ/2)]$. Note that by construction we have the following morphism:
\begin{eqnarray*}
f\colon \cX \too \bbP(\mathrm{Sym}^2(V)) && ([v_1], [v_2])\mapsto [v_1 \otimes v_2 + v_2 \otimes v_1]\,.
\end{eqnarray*}
Given a linear subspace $L\subset \mathrm{Sym}^2(V^\ast)$ of dimension $\leq 3$ when $\mathrm{dim}(V)$ is even, resp. of dimension $\leq 6$ when $\mathrm{dim}(V)$ is odd, consider the associated linear section $\cX_L:=f^{-1}(\bbP(L^\perp))$. Note that such a linear section corresponds to the intersection of $\mathrm{dim}(L)$ bilinear divisors in $\cX$ parametrized by $L$.
\begin{theorem}\label{thm:app2-5}
Assume that $\mathrm{codim}(\cX_L)=\mathrm{dim}(L)$. Under this assumption (which holds for a generic choice of $L$), we have the following implications of conjectures 
\begin{eqnarray*}
\mathrm{C}(Y) \Rightarrow \mathrm{C}_{\mathrm{nc}}(\cX_L) & \text{with} & \mathrm{C}\in \{\mathrm{C}^+, \mathrm{D}, \mathrm{V}, \mathrm{B}, \mathrm{W}, \mathrm{T}^l, \mathrm{T}^p, \mathrm{ST}, \mathrm{K}, \mathrm{S}\}\,,
\end{eqnarray*}
where $Y$ is a certain double cover of $\bbP(L)$ induced by $f$. 
\end{theorem}
Intuitively speaking, Theorem \ref{thm:app2-5} shows that in order to solve the noncommutative counterparts of the celebrated conjectures for intersections of bilinear divisors, it suffices to solve the celebrated conjectures for a certain double cover of the projective space $\bbP(L)$. Theorem \ref{thm:app2-5} was proved in \cite[Thm.~1.13]{NC-zero-CD} in the case of the conjectures $\mathrm{C}^+$ and $\mathrm{D}$. The proof of the other cases is similar.
\begin{corollary}[Intersections of up-to-five bilinear divisors]
Let $\cX_L$ be as in Theorem \ref{thm:app2-5}. When $\mathrm{dim}(L)\leq 2$, the following conjectures hold:
\begin{eqnarray*}
\mathrm{C}_{\mathrm{nc}}(\cX_L) & \mathrm{with} & \mathrm{C} \in \{\mathrm{C}^+, \mathrm{D}, \mathrm{V}, \mathrm{B}, \mathrm{W}, \mathrm{T}^l, \mathrm{T}^p, \mathrm{ST}, \mathrm{P}, \mathrm{K}, \mathrm{S}\} \,.
\end{eqnarray*}
Moreover, $\mathrm{C}_{\mathrm{nc}}(\cX_L)$ holds when $\mathrm{dim}(L)\leq 3$, $\mathrm{D}_{\mathrm{nc}}(\cX_L)$ holds when $\mathrm{dim}(L)\leq 3$ or when $\mathrm{dim}(L)\leq 5$ and $p=0$, and $\mathrm{V}_{\mathrm{nc}}(\cX_L)$ holds when $\mathrm{dim}(L)\leq 3$.
\end{corollary}

%-------------------------------------------------------------------------------
\subsection{Moishezon manifolds associated to quartic double solids}
%-------------------------------------------------------------------------------
Let $k=\bbC$ be the field of complex numbers and $X \to \bbP^2$ one of the quartic double solids introduced by Artin-Mumford in \cite{Artin-Mumford}. These are examples of unirational, but not rational, conic bundles. Thanks to the work of Cossec \cite{Cossec}, these conic bundles can be alternatively described as those singular double coverings $X \to \bbP^3$ which are ramified over a quartic symmetroid $Z$. On the one hand, we can consider the Enriques surface $S_Z$ obtained as the quotient of a natural involution (acting without fixed points) on the blow-up of $Z$. On the other hand, we can consider a small resolution of singularities $\cX \to X$. Such a resolution is not an algebraic variety, but rather a Moishezon manifold\footnote{Recall that a {\em Moishezon manifold} $\cX$ is a compact complex manifold whose field of meromorphic functions on each component has transcendence degree equal to the dimension of the connected component. As proved by Moishezon in \cite{Moishezon}, $\cX$ is a smooth projective $\bbC$-scheme if and only if it admits a K\"ahler metric. In the remaining cases, as proved by Artin in \cite{Artin}, $\cX$ is a proper algebraic space over $\bbC$.}.
\begin{theorem}\label{thm:app2-6}
We have the equivalences of conjectures
\begin{eqnarray*}
\mathrm{C}(S_Z) \Leftrightarrow \mathrm{C}_{\mathrm{nc}}(\cX) &\text{with} & \mathrm{C} \in \{\mathrm{C}^+, \mathrm{D}, \mathrm{V}, \mathrm{S}\}
\end{eqnarray*}
as well as the implication $\mathrm{K}(S_Z) \Rightarrow \mathrm{K}_{\mathrm{nc}}(\cX)$.
\end{theorem}
Roughly speaking, Theorem \ref{thm:app2-6} shows that in what concerns the celebrated conjectures, there is no difference between Enriques surfaces and Moishezon manifolds. Theorem \ref{thm:app2-6} was proved in \cite[Thm.~1.14]{Crelle} in the case of the conjecture $\mathrm{V}$. The proof of the other cases is similar.
\begin{corollary}
The conjectures $\mathrm{C}_{\mathrm{nc}}(\cX)$, with $\mathrm{C} \in \{\mathrm{C}^+, \mathrm{D}, \mathrm{V}\}$, hold.
\end{corollary}


\begin{thebibliography}{00}

\bibitem{Artin} M. Artin, {\em Algebraization of formal moduli, II. Existence of modification}. Ann. of Math. {\bf 91}(2) (1970), 88--135.

\bibitem{Artin-Mumford} M. Artin and D. Mumford, {\em Some elementary examples of unirational varieties which are not rational}. Proc. London Math. Soc. (3) {\bf 25} (1972), 75--95.

\bibitem{Bass} H.~Bass, {\em Some problems in classical algebraic $K$-theory}. Algebraic K-theory, II: ``Classical'' algebraic $K$-theory and connections with arithmetic (Proc. Conf., Battelle Memorial Inst., Seattle, Wash., 1972), pp. 3--73. LNM {\bf 342}, 1973.

\bibitem{Beilinson} A.~Beilinson, {\em Height pairing between algebraic cycles}. K-theory, arithmetic and geometry (Moscow, 1984--1986), 1--25, 
Lecture Notes in Math., {\bf 1289}, Springer, Berlin, 1987.

\bibitem {Crelle} M.~Bernardara, M.~Marcolli, and G.~Tabuada {\em Some remarks concerning Voevodsky's nilpotence conjecture}. Journal f\"ur die reine und angewandte Mathematik, {\bf 738} (2018), 299--312.

\bibitem{Berthelot} P. Berthelot, {\em Cohomologie cristalline des sch\'emas de caract\'eristique $p>0$}. Lecture Notes in Math. {\bf 407}, Springer-Verlag, New York, 1974.

\bibitem{Bloch} S.~Bloch, {\em Lectures on algebraic cycles}. Duke University Mathematics Series, IV. Duke University, Mathematics Department, Durham, N.C., 1980. 

\bibitem{BO} A.~Bondal and D.~Orlov, {\em Semiorthogonal decomposition for algebraic varieties}. Available at arXiv:alg-geom/9506012.

\bibitem{Cadman} C. Cadman, {\em Using stacks to impose tangency conditions on curves}. Amer. J. Math. {\bf 129}
(2007), no.~2, 405--427.

\bibitem{Cossec} F. Cossec, {\em Reye congruences}. Trans. AMS, Vol. {\bf 280} no.~2 (1983), p. 737--751.

\bibitem{Deligne-Moscow} P.~Deligne, {\em Cat\'egories tensorielles}. Mosc. Math. J. {\bf 2} (2002), no. 2, 227--248. Dedicated to Yuri I. Manin on the occasion of his 65th birthday.

\bibitem{Deligne} \bysame, {\em La conjecture de Weil I}. Inst. Hautes \'Etudes Sci. Publ. Math. {\bf 43} (1974), 273--307.

\bibitem{Elment} E.~Elmanto, {\em Topological periodic cyclic homology of smooth $\bbF_p$-algebras}. Talk at the Arbeitsgemeinschaft 2018. Available at \url{https://www.mfo.de/occasion/1814/www_view}.

\bibitem{FT} B.~Feigin and B.~Tsygan, {\em Additive $K$-theory, $K$-theory, arithmetic and geometry (Moscow, 1984--1986)}. Lecture Notes in Math., vol. {\bf 1289}, Springer, Berlin, 1987, pp. 67--209.

\bibitem{Geisser} T.~Geisser, {\em Tate's conjecture, algebraic cycles and rational $K$-theory in characteristic $p$}. $K$-Theory {\bf 13} (1998), no.~2, 109--122.

\bibitem{Quillen} D. Grayson, {\em Finite generation of $K$-groups of a curve over a finite field (after Daniel
Quillen)}. Algebraic $K$-theory, Part I (Oberwolfach, 1980), pp. 69--90, LNM {\bf 966}, 1982.

\bibitem{Grothendieck} A.~Grothendieck, {\em Standard conjectures on algebraic cycles}. 1969 Algebraic Geometry (Internat. Colloq., Tata Inst. Fund. Res., Bombay, 1968) pp. 193--199. Oxford Univ. Press, London.

\bibitem{Grothendieck1} \bysame, {\em Formule de Lefschetz et rationalit\'e des fonctions $L$}. S\'eminaire Bourbaki {\bf 279} (1965).

\bibitem{Guletskii} V.~Guletskii, {\em Finite-dimensional objects in distinguished triangles}. J. Number Theory {\bf 119} (2006), no. 1, 99--127.

\bibitem{GP} V.~Guletskii and C.~Pedrini, {\em Finite-dimensional motives and the conjectures of Beilinson and Murre}. Special issue in honor of Hyman Bass on his seventieth birthday. Part III. $K$-Theory {\bf 30} (2003), no. 3, 243--263. 

\bibitem{Harder} G.~Harder, {\em Die Kohomologie $S$-arithmetischer Gruppen \"uber Funktionenk\"orpern}. Invent.
Math. {\bf 42} (1977), 135--175.

\bibitem{Lars} L. Hesselholt, {\em Topological periodic cyclic homology and the Hasse-Weil zeta function}. Con-
temporary Mathematics {\bf 708} (2018), 157--180.

\bibitem{IU} A. Ishii and K. Ueda, {\em The special McKay correspondence and exceptional collections}. Tohoku
Math. J. (2) {\bf 67}(4) (2015), 585--609.

\bibitem{Kahn1} B.~Kahn, {\em \'Equivalences rationnelle et num\'erique sur certaines vari\'et\'es de type ab\'elien sur un corps fini}. Ann. Sci. Ec. Norm. Sup. {\bf 36} (2003), 977--2002.

\bibitem{KS} B.~Kahn and R.~Sebastian, {\em Smash-nilpotent cycles on abelian $3$-folds}. Math. Res. Lett. {\bf 16} (2009), no. {\bf 6}, 1007--1010.

\bibitem{KM} N.~Katz and W.~Messing, {\em Some consequences of the Riemann hypothesis for varieties over finite fields}. Invent. Math. {\bf 23} (1974), 73--77.

\bibitem{Kedlaya} K.~Kedlaya, {\em Fourier transforms and $p$-adic ``Weil II''}. Compos. Math. {\bf 142} (2006), no.~6, 1426--1450.

\bibitem{Keller} B.~Keller, {\em On differential graded categories}. International Congress of Mathematicians. Vol. II, Eur. Math. Soc., Z\"urich, 2006, pp. 151--190.

\bibitem{Kimura} S.-I.~Kimura, {\em Chow groups are finite dimensional, in some sense}. Math. Ann. {\bf 331} (2005), no.~1, 173--201.

\bibitem{Kleiman} S.~L.~Kleiman, {\em The standard conjectures}. Motives (Seattle, WA, 1991), 3--20, Proc. Sympos. Pure Math., {\bf 55}, Part 1, Amer. Math. Soc., Providence, RI, 1994.

\bibitem{Kleiman1} \bysame, {\em Algebraic cycles and the Weil conjectures}. Dix expos\'es sur la cohomologie des
sch\'emas, 359--386, Adv. Stud. Pure Math., {\bf 3}, North-Holland, Amsterdam, 1968.

\bibitem{Miami} Maxim Kontsevich, {\em Mixed noncommutative motives}. Talk at the Workshop on Homological Mirror Symmetry, Miami, 2010. Available at \url{www-math.mit.edu/auroux/frg/miami10-notes}.

\bibitem{finMot} \bysame, {\em Notes on motives in finite characteristic}. Algebra, arithmetic, and geometry: in honor of Yu. I. Manin. Vol. II, Progr. Math., vol. {\bf 270}, Birkh\"auser Boston, Inc., Boston, MA, 2009, pp. 213--247.

\bibitem{Lefschetz} \bysame, {\em Hodge structures in non-commutative geometry}. XI Solomon Lefschetz Memorial Lecture series. Contemp. Math., {\bf 462}, Non-commutative geometry in mathematics and physics, 1--21, Amer. Math. Soc., Providence, RI, 2008. 

\bibitem{Kontsevich-talk} \bysame, {\em Categorification, NC Motives, Geometric Langlands, and Lattice Models}. Talk at the Geometric Langlands Seminar, University of Chicago, 2006. Notes available at the webpage \url{https://www.ma.utexas.edu/users/benzvi/notes.html}.

\bibitem{IAS} \bysame, {\em Noncommutative motives}. Talk at the IAS on the occasion of the $61^{\mathrm{st}}$ birthday of Pierre Deligne (2005). Available at \url{http://video.ias.edu/Geometry-and-Arithmetic}.

\bibitem{Kuznetsov-CY} A.~Kuznetsov, {\em Calabi-Yau and fractional Calabi-Yau categories}. Available at arXiv:1509.07657. To
appear in J. Reine Angew. Math.

\bibitem{delPezzo} A. Kuznetsov, {\em Derived categories of families of sextic del Pezzo surfaces}. Available at
arXiv:1708.00522.

\bibitem{Lieberman} D.~Lieberman, {\em Numerical and homological equivalence of algebraic cycles on Hodge manifolds}. Amer. J. Math. {\bf 90}, 366--374, 1968.

\bibitem{Manin} Y.~Manin, {\em Correspondences, motifs and monoidal transformations}. Mat. Sb. (N.S.) {\bf 77} (119) (1968), 475--507.

\bibitem{JEMS} M.~Marcolli and G.~Tabuada, {\em Noncommutative numerical motives, Tannakian structures, and motivic Galois groups}. Journal of the European Mathematical Society {\bf 18} (2016), 623--655.

\bibitem{Mazza} C.~Mazza, {\em Schur functors and motives}. $K$-Theory {\bf 33} (2004), no.~2, 89--106.

\bibitem{Milne} J.~Milne, {\em The Tate conjecture over finite fields}. AIM talk. Available at Milne's webpage \url{http://www.jmilne.org/math/articles/2007e.pdf.}

\bibitem{Moishezon} B.~G.~Moishezon, {\em On $n$-dimensional compact varieties with $n$ algebraically independent meromorphic functions, I, II and III}. Izv. Akad. Nauk SSSR Ser. Mat., {\bf 30}: 133--174.

\bibitem{Scholze} T. Nikolaus and P. Scholze, {\em On topological cyclic homology}. Available at arXiv:1707.01799. To appear in Acta Math.

\bibitem{Shermenev} A.~Shermenev, {\em The motive of an abelian variety}. Funct. Anal. {\bf 8} (1974), 47--53.

\bibitem{survey} G.~Tabuada, {\em Recent developments on noncommutative motives}. New directions in homotopy theory, 143--173, 
Contemp. Math., {\bf 707}, Amer. Math. Soc., Providence, RI, 2018.

\bibitem{MRL} \bysame, {\em A note on the Schur-finiteness of linear sections}. Mathematical Research Letters, {\bf 25} (2018), no.~1, 237--253.

\bibitem{NC-zero-CD} \bysame, {\em A note on Grothendieck's standard conjectures of type $\mathrm{C}^+$ and $\mathrm{D}$}. Proc. Amer. Math. Soc. {\bf 146} (2018), no.~4, 1389--1399. 

\bibitem{Kimura-curve} \bysame, {\em  Kimura-finiteness of quadric fibrations over smooth curves}. Presented by Christophe Soul\'e. C. R. Math. Acad. Sci. Paris {\bf 355} (2017), no. 6, 628--632.

\bibitem{book} \bysame, {\em Noncommutative motives}. With a preface by Yuri I. Manin. University Lecture Series {\bf 63}. American Mathematical Society, Providence, RI, 2015.

\bibitem{JNCG} \bysame, {\em Chow motives versus noncommutative motives}. J. Noncommut. Geom. {\bf 7} (2013), no. 3, 767--786.

\bibitem{NCWeil} \bysame, {\em Noncommutative Weil conjecture}. Available at arXiv:1808.00950.

\bibitem{NCTate} \bysame, {\em HPD-invariance of the Tate, Beilinson and Parshin conjectures}. Available at arXiv:1712.05397.

\bibitem{NCSchur} \bysame, {\em Schur-finiteness (and Bass-finiteness) conjecture for quadric fibrations and for families of sextic du Val del Pezzo surfaces}. Available at arXiv:1708.05382.

\bibitem{NCpositive-CD} \bysame, {\em On Grothendieck's standard conjectures of type $\mathrm{C}^+$ and $\mathrm{D}$ in positive characteristic}. Available at arXiv:1710.04644.

\bibitem{NCpositive} \bysame, {\em Noncommutative motives in positive characteristic and their applications}. Available at arXiv:1707.04248. To appear in Advances in Mathematics.

\bibitem{Orbifold} G.~Tabuada and M. Van den Bergh, {\em Additive invariants of orbifolds}. Geom. Topol. {\bf 22} (2018), no.~5, 3003--3048. 

\bibitem{Tate1} J.~Tate, {\em Conjectures on algebraic cycles in $l$-adic cohomology}. Motives (Seattle, WA, 1991), 71--83, Proc. Sympos. Pure Math., {\bf 55}, Part 1, Amer. Math. Soc., Providence, RI, 1994. 

\bibitem{Tate} \bysame, {\em Algebraic cycles and poles of zeta functions}. Arithmetical Algebraic Geometry (Proc. Conf. Purdue Univ., 1963) pp. 93--110. Harper \& Row, New York 1965.

\bibitem{Thomason-Tate} R. Thomason, {\em A finiteness condition equivalent to the Tate conjecture over $\bbF_q$}. Algebraic $K$-theory and algebraic number theory (Honolulu, HI, 1987), 385--392, Contemp. Math., {\bf 83}, Amer. Math. Soc., Providence, RI, 1989.

\bibitem{Totaro} B.~Totaro, {\em Recent progress on the Tate conjecture}. Bull. Amer. Math. Soc. {\bf 54} (2017), 575--590.

\bibitem{Voevodsky2} V.~Voevodsky, {\em Triangulated categories of motives over a field}. Cycles, transfers, and motivic homology theories, Ann. of Math. Stud., vol. {\bf 143}, Princeton Univ. Press, 2000, pp. 188--238.

\bibitem{Voevodsky} \bysame, {\em A nilpotence theorem for cycles algebraically equivalent to zero}. Int. Math. Res. Not. IMRN 1995 (1995), no. {\bf 4}, 187--198.

\bibitem{Voisin} C.~Voisin, {\em Remarks on zero-cycles of self-products of varieties}. in: Moduli of vector bundles (Sanda 1994,
 Kyoto 1994), Lecture Notes in Pure and Appl. Math. {\bf 179}, Dekker, New York (1996), 265--285.

\bibitem{Weil} A.~Weil, {\em Vari\'et\'es ab\'eliennes et courbes alg\'ebriques}. Hermann, Paris (1948).

\bibitem{Weil1} \bysame, {\em Sur les courbes alg\'ebriques et les vari\'et\'es qui s'en d\'eduisent}. Actualit\'es Sci. Ind., no. {\bf 1041} = Publ. Inst. Math. Univ. Strasbourg {\bf 7} (1945). Hermann et Cie., Paris, 1948.

\end{thebibliography}
\end{document}